\documentclass[12pt]{article}
\usepackage{amsmath,amssymb,amsthm,setspace,latexcad}
\usepackage{graphicx}
 \newtheorem{theorem}{Theorem}
 
 \newtheorem{lemma}[theorem]{Lemma}

 {\theoremstyle{definition}
 \newtheorem{definition}{Definition}

 }
\def\N{\ensuremath{\mathbb N}} 
\def\R{\ensuremath{\mathbb R}} 
\def\Z{\ensuremath{\mathbb Z}} 

\renewcommand{\bf}{\bfseries}

\newcommand{\dis}{\displaystyle}

\newcommand{\DDD}{\Delta}

\newcommand{\Ics}[1]{I^{**}_{#1}}

\newcommand{\aaa}{\alpha}

\newcommand{\begall}{\begin{align*}}
\newcommand{\begal}{\begin{align}}
\newcommand{\begdef}{\begin{definition}}
\newcommand{\begeqq}{\begin{equation}}

\newcommand{\cai}{{\cal I}}

\newcommand{\cak}{{\cal K}}

\newcommand{\capp}{{\cal P}}

\newcommand{\cat}{{\cal T}}

\newcommand{\cax}{{\cal X}}

\newcommand{\defeq}{{\buildrel {\rm def}\over =}}

\newcommand{\eee}{\epsilon}

\newcommand{\ess}{\emptyset}

\newcommand{\fff}{{\phi}}

\newcommand{\isss}{\hbox{ispt}}
\newcommand{\kkks}{s}

\newcommand{\lf}{\lfloor}
\newcommand{\lxto}{l_{x,t_{1}}}

\newcommand{\nk}{n_{\kkks}}

\newcommand{\ocax}{\overline{{\cal X}}}

\newcommand{\ointt}{\overline{\int}}
\newcommand{\olll}{\overline{\lambda}}
\newcommand{\ooo}{\omega}
\newcommand{\oot}[1]{\left (\frac{1}2\right)^{#1}}
\newcommand{\oo}{\infty}

\newcommand{\ppp}{\psi}

\newcommand{\restr}[1]{\big|_{#1}}

\newcommand{\rf}{\rfloor}
\newcommand{\rrr}{\rho}
\newcommand{\sj}[1]{_{j,#1}}
\newcommand{\sm}{\setminus}

\newcommand{\so}[1]{_{1,#1}}
\newcommand{\spt}{\text{spt}\,}
\newcommand{\sse}{\subset}

\newcommand{\stj}[1]{_{t,j,#1}}
\newcommand{\sto}[1]{_{t_{0},#1}}
\newcommand{\supp}{\text{spt}\,}

\newcommand{\tcs}[1]{t^{*}_{#1}}

\newcommand{\ttt}{\theta}

\renewcommand{\lll}{\lambda}
\title{The $(L^{1},L^{1})$ bilinear Hardy-Littlewood function and Furstenberg averages}
\author{
Idris Assani\thanks{The first author acknowledges support by NSF grant DMS 0456627},
 Department of Mathematics,\\
 University of North Carolina at Chapel Hill,\\
 Chapel Hill, North Carolina 27599, USA\\
email: assani@email.unc.edu\\
{\tt www.math.unc.edu/Faculty/assani}
\\
\\
and  \\
\\
Zolt\'an Buczolich\thanks{
Research supported by the Hungarian
National Foundation for Scientific research T049727.
In December of 2007 this author received partial support from the Erwin Schr\"ordinger
Institute in Vienna to discuss some of these results and to do some work related to this paper.
\newline\indent {\it ESI preprint number 2015}
},
Department of Analysis, E\"otv\"os Lor\'and\\
University, P\'azm\'any P\'eter S\'et\'any 1/c, 1117 Budapest, Hungary\\
email: buczo@cs.elte.hu\\
{\tt www.cs.elte.hu/\hbox{$\sim$}buczo} }

\date{\today}
\begin{document}
\maketitle

\medskip

\begin{abstract}
Let $(X,\mathcal{B}, \mu, T)$ be an ergodic dynamical system on a
 non-atomic finite measure space.
 Consider the maximal function
 $\dis R^*:(f, g) \in L^1\times L^1 \rightarrow R^*(f, g)(x) =
 \sup_{n} \frac{f(T^nx)g(T^{2n}x)}{n}.$
  We show that there exist $f$ and $g$ such that
 $R^*(f, g)(x)$ is not finite almost everywhere. Two consequences
 are derived.
 The bilinear Hardy--Littlewood maximal function fails to be a.e.
 finite for all functions $(f, g)\in L^1\times L^1.$
 The Furstenberg averages do not converge
 for all pairs of $(L^{1},L^{1})$ functions, while by a result
 of J. Bourgain these averages converge for all pairs of
 $(L^{p},L^{q})$ functions with $\frac{1}{p}+\frac{1}{q}\leq 1.$
\end{abstract}

\section{Introduction}
The bilinear Hardy--Littlewood maximal function was introduced by
Alberto Calder\'on in the 1960's. It is defined for $f, g$
measurable functions as
 $$M^*(f, g)(x) = \sup_t\frac{1}{2t}\int_{-t}^t f(x + s) g(x +
 2s)ds.$$ Our purpose is to prove that $M^*$ is not always a.e finite when
 the functions $f$ and $g$ are in $L^1.$
 \begin{theorem} \label{ith0}
 There exist functions $f, g$ both in $L^1$ for which the bilinear
 Hardy-Littlewood maximal function
 $$M^*(f, g)(x) = \sup_t\frac{1}{2t}\int_{-t}^t f(x + s) g(x +
 2s)ds$$
 is not a.e. finite.
 \end{theorem}
 To prove this theorem we use Ergodic Theory.
  The Ergodic Theory version of the bilinear Hardy--Littlewood
maximal function is defined for $f,g\geq 0$ as
$$
 \mathcal{M}(f,g)(x)= \sup_N\frac{1}{2N+1}\sum_{n=-N}^N f(T^nx)g(T^{2n}x)
  ,$$
where $T$ is an ergodic measure preserving transformation of a
non-atomic probability measure space, $f\in L^{p}$, and $g\in
L^{q}$. They are called Furstenberg averages. They appear in H.
Furstenberg's paper in 1977 \cite{Furstenberg}.

A transference argument shows that the class of functions in
$L^p\times L^q$, for which $M^*$ and $\mathcal{M}$ are a.e. finite
is the same. For $f, g\geq 0$ we have
$$\mathcal{M}(f,g)(x)\geq  \sup_N\frac{ f(T^Nx)g(T^{2N}x)}{2N+1},$$
  the tail of the averages  $\mathcal{M}(f,g)(x).$ We consider the
maximal function
 $$R^*(f, g)(x) = \sup_n\frac{f(T^nx)g(T^{2n}(x))}{n}.$$

In \cite{[fgtb]} we showed that for all $p, q \geq 1$ such that $\frac{1}{p} +
 \frac{1}{q} <2,$ $R^*$ maps $L^p\times L^q$ into $L^r$ as soon
 as $0<r<1/2.$ This implies that $R^{*}(f,g)$ is finite almost everywhere
 and $\frac{f(T^nx)g(T^{2n}x)}{n}\to 0$ for a.e. $x$ as $n\to\oo$.

In this paper we show that for $p=q=1,$  $R^*(f,g)(x)$ is not
 finite almost everywhere for all $f$ and $g$.

 \begin{theorem}\label{ith1}
 Let $(X, \mathcal{B}, \mu, T)$ be an ergodic measure
 preserving transformation on a finite non-atomic measure
 space. Then there exist functions $f, g$ both in $L^1_+(X)$ for
 which the maximal function $$R^*(f, g)(x) = \sup_n\frac{f(T^nx)g(T^{2n}x)}{n}$$
 is not finite a.e.
 \end{theorem}

  The example of the identity map shows that Theorem \ref{ith1} is false
without the ergodicity assumption.

 Theorem \ref{ith1} gives us three conclusions.
 First it solves an open problem in Ergodic Theory.
Indeed, a deep result of J. Bourgain, \cite{Bourgain}, showed that
the Furstenberg averages converge a.e. as soon as the H\"olderian
duality is respected, (i.e. $\frac{1}{p}+\frac{1}{q} \leq 1$).
  Theorem \ref{ith1} shows that
 these averages do not converge for all pairs of $(L^1, L^1)$
 functions as the tail of these averages does not converge a.e.
 to zero for some functions $f,g\in L^1$. This is the content of the
 following result
 \begin{theorem} \label{ith0'}Given an ergodic measure preserving transformation $T,$ on a
 nonatomic probabilty measure space,  we can find functions $f,g\in
 L^1$ for which the Furstenberg averages
  $$\frac{1}{N}\sum_{n=0}^N f(T^nx)g(T^{2n}x)
  $$ do not converge a.e.
 \end{theorem}
 Secondly, by transference the unboundedness of $\mathcal{M}(f,g)(x)$ implies
 the same result for the bilinear Hardy Littlewood maximal function
 in $(L^1\times L^1)$ and gives a proof of Theorem \ref{ith0}.

 A third consequence of Theorem \ref{ith1} is that 1/2 is an optimal bound for $R^*.$
\vskip1ex
 In view of all these three consequences we just need to focus on proving Theorem \ref{ith1}.
 \vskip1ex

 Let us fix some
notation. Given $\fff:\R\to\R$, periodic by $p$ we put
$$\ointt \fff=\frac{1}p \int_{0}^{p}\fff(x)dx.$$
Given a Lebesgue measurable set $A$, periodic by
$p$ we put $$\olll(A)=\frac{1}p\lambda (A\cap [0,p)).$$

For a function $\fff:\R\to \R$, we will denote by $\spt(\fff)$
the set of those $x$'s for which
$\fff(x)\not=0$. This notation differs slightly from the support of a function.
The functions for which we will apply it will be constant
on intervals of the form $[k,k+1)$, $k\in \Z$ and hence
$\spt \fff$ differs only by some endpoints of these intervals
from the set which is usually considered to be the support
of a function.

\section{Main Results}

We want to prove Theorem \ref{ith1}. To this end we introduce the following
definition.

Suppose that $\capp \sse \N$ is an infinite set, $\frac{1}{2}<\rho
<1,$ $0<\epsilon <\epsilon _{0}\defeq 1/10$, $\kkks,\, \aaa\in \N.$
We also suppose that $M\geq 2, M\in \N$

We will define the functions $\fff_{i},\ \ppp_{k}:\R\to\R$ so that
if $\DDD$ is any of these functions then $\DDD(x)=\DDD(\lf x \rf)$.
\vskip1ex

 We fix the parameters  $\rho,$ $\epsilon$ and $M.$

\begin{definition}\label{deff1}
An $\aaa-\capp-\kkks$-family with $\ppp$-interval $[\aaa,\ooo]$
consists of an integer $n_{\kkks}$ and of non-negative functions
$\fff_{i}$, $i=1,...,n_{\kkks}$, $\ppp_{k}$, $k=\aaa,...,\ooo$ which
are periodic by $p_{\kkks}\in\capp$, $p_{\kkks}\geq\ooo$ and
\begin{align}\label{*f1}
\olll \bigg  \{ x: \max_{k\in[\aaa,\ooo]}\max_{l\leq p_{\kkks}}
&
\frac{\sum_{i=1}^{n_{\kkks}}\fff_{i}(x+l)\ppp_{k}(x+2l)}{l}\geq 1
 \bigg \} >\\ &\nonumber \min \bigg
\{ \eee,\frac{\kkks(M-1)\eee2^{-M}}{2048}
 \bigg
 \},
\end{align}
\begin{equation}\label{**f1}
\kkks\cdot 2^{-M-1}<\sum_{i=1}^{n_{\kkks}}\ointt\fff_{i}<
\kkks\cdot 2^{-M+1},
\end{equation}
\begin{equation}\label{***f1}
\rrr\eee<\sum_{k=\aaa}^{\ooo}\ointt\ppp_{k}<\eee.
\end{equation}
\end{definition}
  Note that in this definition $\aaa$ is predetermined but $\ooo$ is not.
  The number $\ooo$ depends on $\aaa$ and on other objects
of our construction like the functions $\fff$ and  $\ppp.$

  The next lemma shows that the notion of $\aaa-\capp-\kkks$-family
  is ``generic" among infinite sets $\capp$ of positive integers in
  the following sense: if we can find one $\aaa-\capp'-\kkks$-family
  then for any other infinite subset $\capp$ we will be able to find
  an $\aaa-\capp-\kkks$-family.
  This lemma allows us to restrict the proof of the existence of such families to infinite set
  of integers consisting only of powers of two.

\begin{lemma}\label{lemf3}
Suppose that $\capp'\sse \N$ is an arbitrary infinite set and there
exists an $\aaa-\capp'-\kkks$-family with $\ppp'$-interval
$[\aaa,\ooo].$ Then for an arbitrary other infinite set
$\capp\sse\N$ there exists an $\aaa-\capp-\kkks$-family with
$\ppp$-interval $[\aaa,\ooo]$ as well. Moreover,
\begin{equation}\label{*5*5}
\ointt\fff_{i}'\approx \ointt\fff_{i}, \text{ in fact }\left (1-\frac{p_{s}'}{p_{s}} \right )
\ointt \fff_{i}'\leq \ointt \fff_{i}\leq \ointt\fff_{i}',
\end{equation}
and
\begin{equation}\label{*f16}
\ointt\ppp_{k}'\approx \ointt\ppp_{k}, \text{ in fact }\left (1-\frac{p_{s}'}{p_{s}} \right )
\ointt \ppp_{k}'\leq \ointt \ppp_{k}\leq \ointt\ppp_{k}'
\end{equation}
where $\fff_{i}'$, $\ppp_{k}'$ belongs to the ``old" and $\fff_{i}$, $\ppp_{k}$
to the ``new" family and these families are periodic by $p_{s}'$ and $p_{s}$, respectively.
\end{lemma}

\begin{proof}
Suppose we have $\fff_{i}',$ $i=1,...,n_{\kkks}$,
$\ppp_{k}'$, $k=\aaa,...,\ooo$ periodic by $p_{\kkks}'
\in\capp'$, $p_{\kkks}'\geq\ooo$
satisfying (\ref{*f1}-\ref{***f1}).
Since $\capp$ contains infinitely many terms there
are arbitrarily large elements $p_{\kkks}\in \capp$.
We will select a sufficiently large $p_{\kkks}\gg
p_{\kkks}'.$

To define $\fff_{i}$, $\ppp_{k}$ periodic by $p_{\kkks}$
it is sufficient to define them on $[0,p_{\kkks})$
and then extend their definition onto $\R$ by
periodicity.
If $x\in [0, \lf p_{\kkks}/p_{\kkks}' \rf\cdot p_{\kkks}')$
then set $\fff_{i}(x)=\fff_{i}'(x),$ $i=1,...,n_{\kkks}$
and $\ppp_{k}(x)=\ppp'_{k}(x),$ $k=\aaa,...,\ooo.$
If $x\in [\lf p_{\kkks}/p_{\kkks}' \rf\cdot p_{\kkks}',
p_{\kkks})$ then set $\fff_{i}(x)=0=\ppp_{k}(x)$
for $i=1,...,n_{\kkks}$ and $k=\aaa,...,\ooo.$
Since in (\ref{*f1}-\ref{***f1}) concerning
$\olll$ and $\ointt$ there are strict inequalities
it is not difficult to see that (\ref{*f1}-\ref{***f1})
hold for $\fff_{i}$ and $\ppp_{k}$ if $p_{\kkks}$
is sufficiently large.  For example, we show that \eqref{***f1}
and \eqref{*f16} hold.
We have
\begin{equation}\label{*45esta}
\ointt\ppp_{k}=\frac{1}{p_{s}}\int_{[0,p_{s})}\ppp_{k}
=\frac{1}{p_{s}}\int_{[0,\lf p_{s}/p_{s}' \rf p_{s}')}\ppp_{k}'=
\end{equation}
$$\frac{\lf p_{s}/p_{s}' \rf p_{s}'}{p_{s}}\frac{1}{\lf p_{s}/p_{s}' \rf p_{s}'}
\int_{[0,\lf p_{s}/p_{s}' \rf p_{s}')}\ppp_{k}'=
\frac{\lf p_{s}/p_{s}' \rf p_{s}'}{p_{s}}\ointt\ppp_{k}'.
$$
From this and
\begin{equation}\label{*45estb}
\frac{p_{s}-p_{s}'}{p_{s}}<\frac{\lf p_{s}/p_{s}' \rf p_{s}'}{p_{s}}\leq 1
\end{equation}
it follows \eqref{*f16} and adding these inequalities for $k=\aaa$
to $\ooo$ we obtain \eqref{***f1} when $p_{s}$ is sufficiently large.
\end{proof}

The next simple ``independence lemma" will be useful later.

\begin{lemma}\label{*indep}
Suppose $t,n, \pi_{1},\pi_{2}, \ttt_{1},\ttt_{2}\in \Z$,
$\ttt_{1},\ttt_{2}\geq 0.$ Consider two sets  $\cax_{1},\cax_{2}\sse
[t2^{n},(t+1)2^{n})$ with the following
properties: \\
 a) They are ``periodic within"
$[t2^{n}+\ttt_{1},(t+1)2^{n}-\ttt_{2})$ by $\pi_{1}$ and $\pi_{2}$,
respectively. This means that for $i=1,2$ if $x,x+\pi_{i}\in
[t2^{n}+\ttt_{1},(t+1)2^{n}-\ttt_{2})$
then $x\in \cax_{i}$ iff $x+\pi_{i}\in \cax_{i}$.\\
b) They consist of integer intervals. This means that $x\in\cax_{i}$
iff $\lf x\rf\in \cax_{i}$
for $i=1,2.$\\
Then for relatively prime $\pi_{1}$ and $\pi_{2}$ we have
\begin{equation}\label{*0pot}
\lll(\cax_{1}\cap\cax_{2})<2\frac{\lll(\cax_{1})\lll(\cax_{2})}{2^{n}}
\end{equation}
if
 $2^{n}$ is much larger than $\max \{ \pi_{1},\pi_{2},\ttt_{1},
\ttt_{2}  \}$.
\end{lemma}

\begin{proof}
Set $X_{1}=\cax_{1}\cap \Z$, $X_{2}=\cax_{2}\cap \Z$. Since the sets
$\cax_{1}$ and $\cax_{2}$ consist of integer intervals
$\lll(\cax_{i})=\# X_{i}$, ($i=1,2$) and $\lll(\cax_{1}\cap
\cax_{2})=\# (X_{1}\cap X_{2}).$

Suppose $i\in \{ 1,2  \}$ and
\begin{equation}\label{*6bb}
t2^{n}+\ttt_{1}\leq k\pi_{i}<(k+1)\pi_{i}-1< (t+1)2^{n}-\ttt_{2}.
\end{equation}
Then $$L_{i}\defeq \frac{\#(X_{i}\cap
[k\pi_{i},(k+1)\pi_{i}))}{\pi_{i}}$$ does not depend on the choice of
$k$
as long as \eqref{*6bb} is satisfied.
 If $n\to\oo$ with $\pi_{1},$ $\pi_{2}$, $\ttt_{1}$ and $\ttt_{2}$ fixed
 then
$$\frac{\lll(\cax_{i})}{2^{n}}
= \frac{\# X_{i}}{2^{n}}\to L_{i}.$$ Therefore, if
$2^{n}$ is much larger than $\pi_{1},$ $\pi_{2}$, $\ttt_{1}$ and
$\ttt_{2}$ then
\begin{equation}\label{*1pot}
L_{i}<\root 3 \of 2 \frac{\lll(\cax_{i})}{2^{n}} \text{ holds for }i=1,2.
\end{equation}

Suppose
$$t2^{n}+\ttt_{1}\leq k\pi_{1}
\pi_{2}<(k+1)\pi_{1}\pi_{2}-1< (t+1)2^{n}-\ttt_{2}.$$ Since
$\pi_{1}$ and $\pi_{2}$ are relatively prime
$$L_{1}L_{2}=
\frac{\#(X_{1}\cap X_{2}\cap
[k\pi_{1}\pi_{2},(k+1)\pi_{1}\pi_{2}))}{\pi_{1} \pi_{2}}.$$ If
$n\to\oo$ with $\pi_{1},$ $\pi_{2}$, $\ttt_{1}$ and $\ttt_{2}$ fixed
then
$$\frac{\lll(\cax_{1}\cap \cax_{2})}{2^{n}}
= \frac{\# (X_{1}\cap X_{2})}{2^{n}} \to L_{1}L_{2}.$$ Therefore, if
$2^{n}$ is much larger than $\pi_{1},$ $\pi_{2}$, $\ttt_{1}$ and
$\ttt_{2}$ then
$$
 \frac{\lll(\cax_{1}\cap \cax_{2})}{2^{n}}< L_{1}L_{2} \root 3 \of 2.
$$

By using \eqref{*1pot} we have
$$\frac{\lll(\cax_{1}\cap\cax_{2})}{2^{n}}<L_{1}L_{2}
\root 3 \of 2<(\root 3 \of 2\frac{\lll(\cax_{1})}{2^{n}}) (\root 3
\of 2\frac{\lll(\cax_{2})}{2^{n}}) \root 3 \of 2
$$
and after multiplying by $2^{n}$ this implies \eqref{*0pot}.

\end{proof}
   Theorem \ref{ith1} is a consequence of Theorem \ref{f2}. The proof is quite
   long. Therefore we have divided it into two steps
and several substeps (14 SUBSTEPS).
   STEP 1 contains six substeps that correspond to the first
   induction step for $s=1.$ In STEP 2 we finish the induction
   argument on $s.$
In Section \ref{derith1} we show how this result allows
   us to derive a proof of Theorem \ref{ith1}.
We recall that the parameters $\frac{1}{2}<\rho<1,$
$0<\epsilon<\frac{1}{10}=\epsilon_0,$ and $M\in \N, M\geq 2,$ are
fixed in the definition of an $\aaa-\capp-\kkks$-family.

\begin{theorem}\label{f2}
Suppose that $\capp\sse\N$ is an infinite set. Then for every
$\aaa,$ $s\in \N$ there exists an $\aaa-\capp-\kkks$-family.
\end{theorem}

\begin{proof}
We do mathematical induction on $\kkks.$
By Lemma \ref{lemf3} we can suppose that $\capp$
consists of powers of $2$.

\subsection{STEP 1: The $s=1$ case\\
SUBSTEP 1a: Interval supports}\label{subsubstep1a}

First we show that for any possible choice of $\aaa$, $\capp,$
$\frac{1}{2}<\rrr<1$, $0<\eee<\frac{1}{10}$, and $M\in \N, M\geq 2,$
one can find $\aaa-\capp-1$-families.

We will select later a suitably large $n_{1}$ and for $i<n_{1}$ we
set $\fff_{i}(x)=0$ for all $x\in\R$. To define $\fff_{n_{1}}$ and
the functions $\ppp_{k}$, $k=\aaa,...,\ooo$ one could come up with a
somewhat simpler definition but to help the reading of the more
technical later steps of our induction we introduce already at the
first step some of the features used later.

We choose integers $0=k_{0}<k_{1}<...<k_{M}$
so that $k_{1}\gg \aaa$.

The interval support of $\fff_{n_{1}}$ at level $k$
is defined as
\begin{align}\label{*iss*}
\isss_{k}(\fff_{n_{1}})=  \bigcup  &\{ [(t-1)2^{n_{1}+k},t2^{n_{1}+k})
:\\
\nonumber
&t\in \Z,\   \spt(\fff_{n_{1}})\cap  [(t-1)2^{n_{1}+k},t2^{n_{1}+k})
\not=\ess \} .
\end{align}

\placedrawing[h]{fgt1.lp}{One component of $\isss_{k_{j}}(\fff_{n_{1}})$}{figure1}

Outside its support $\fff_{n_{1}} $ will vanish, so by giving
its interval support at different levels we can define it.
We zoom in during this definition. (This is similar
to the procedure of defining the triadic Cantor set
as the intersection of closed sets at level $k$ consisting
of $2^{k}$ many intervals of length $3^{-k}$. Though in our
construction we will use only finitely many steps of
zooming in.)
We set $$\isss_{k_{M}}(\fff_{n_{1}} )=\R.$$
\begin{align}\nonumber
&\text{For }j\in  \{ 1,...,M  \} \text{ an interval }
[(t-1)2^{n_{1}+k_{j-1}},t 2^{n_{1}+k_{j-1}})\sse
\isss_{k_{j}}(\fff_{n_{1}} )
 \\
&\label{*isptdef}\text{belongs to }
\isss_{k_{j-1}}(\fff_{n_{1}} ) \text{ if and only if }
t \text{ is even.}\end{align}
(See Figure \ref{figure1}.)
This implies
\begin{equation}\label{*ISP0*}
\olll(\isss_{k_{j}}(\fff_{n_{1}}))=\left (\frac{1}{2} \right )^{M-j},
\end{equation}
and
\begin{equation}\label{*ISP*}
\olll(\isss_{k_{j-1}}(\fff_{n_{1}}))=\olll(\isss_{k_{j}}(\fff_{n_{1}})\sm
\isss_{k_{j-1}}(\fff_{n_{1}}))=\frac{1}{2}\olll(\isss_{k_{j}}(\fff_{n_{1}})).
\end{equation}

 It will be useful to keep in mind
for further reference that $t-1$, which corresponds
to the left endpoint of a support interval is odd.

\subsection{SUBSTEP 1b: Definitions of
$\fff_{n_{1}}$ and $\ppp_{k_{j}'}$}

We will define $\fff_{n_{1}} $ so that it will be periodic by
$2^{n_{1}+k_{M}}\in \capp$.

If $[(t-1)2^{n_{1}},t 2^{n_{1}})\sse\isss_{0}(\fff_{n_{1}})=\isss_{k_{0}}(\fff_{n_{1}})$
then
\begin{align}\label{*f5}
&\fff_{n_{1}}(x)=2^{n_{1}}\text{ if }x\in [(t-1)2^{n_{1}},(t-1)2^{n_{1}}
+1),\text{ and }\\
&\fff_{n_{1}}(x)=0\text{ if }x\in [(t-1)2^{n_{1}}+1,t\cdot 2^{n_{1}}
).
\end{align}

From \eqref{*ISP0*} and \eqref{*ISP*} used with $j=1$ it follows that
\begin{equation}\label{*AVINTPHIN1*}
\ointt\sum_{i=1}^{n_{1}}\fff_{i}=\ointt\fff_{n_{1}}=2^{-M}.
\end{equation}

Set $k_{1}'=\aaa.$

For any $k\not=k_{1}'$, $k\in \N$ and $x\in\isss_{k_{1}}(\fff_{n_{1}})$
we set $\ppp_{k}(x)=0.$

For $x\not\in \isss_{k_{1}}(\fff_{n_{1}})$ we set
$\ppp_{k_{1}'}(x)=0$, that is,
\begin{equation}\label{*psiksup1*}
\ppp_{k_{1}'} \text{ is supported in } \isss_{k_{1}}(\fff_{n_{1}}).
\end{equation}

We will choose $\ppp_{k_{1}'}$ so that it is constant
on intervals of the form $[n,n+1)$ for all $n\in\Z$,
it is periodic by $2^{n_{1}+k_{M}}$,
its range is $ \{ 2^{k_{1}'},0  \}$ and
\begin{equation}\label{2*f6b}
\frac{1+\rrr}{2}\eee\cdot 2^{n_{1}}<
\int_{[(t-1)2^{n_{1}}+1,t2^{n_{1}}+1) }
\ppp_{k_{1}'}(x)dx <\eee 2^{n_{1}},
\end{equation}
provided $[(t-1)2^{n_{1}},t2^{n_{1}})\sse
\isss_{k_{1}}(\fff_{n_{1}})$.
Inequality \eqref{2*f6b}
can be achieved if $2^{n_{1}}$ is sufficiently large.

Suppose $j\geq 2$ and we can also suppose that
$\aaa\ll k_{1}\ll ... \ll k_{M}$.

We will choose $n_{1}$ so that $n_{1}\gg k_{M}$ and
$k_{j}'$ such that $k_{j-1}<  k_{j}'\defeq k_{j-1}+10 \ll k_{j}.$
Our constants will be selected by induction in the following order:
$k_{1}$, $k_{1}'$, $k_{2}$, ... $k_{j-1},$ $k_{j}'$, $k_{j}$, ..., $k_{M}$,
$n_{1}$. Assumptions about $k_{j}$'s and $k_{j}'$'s are given in
Section \ref{213}, we emphasize that in \eqref{2*f7a} the value of the fraction
does not depend on $n_{1}$. While we make assumptions about $n_{1}$
at \eqref{2*f6b}, \eqref{*g7}, \eqref{2*f6} and \eqref{*f7g}.

The function
\begin{equation}\label{*psiksup}
\ppp_{k_{j}'} \text{ will be supported in }
\isss_{k_{j}}(\fff_{n_{1}})\sm \isss_{k_{j-1}}(\fff_{n_{1}}).
\end{equation}
At the $s=1$ case one nonzero $\ppp_{k_{j}'}$ is sufficient,
that is, we set
$\ppp_{k}(x)=0$
for all $x\in
\isss_{k_{j}}(\fff_{n_{1}})\sm \isss_{k_{j-1}}(\fff_{n_{1}})$,
 $k\in\N$, $k\not=k_{j}'$.

To each $t\in \Z$
we associate a period $\pi(t)\gg 2^{k_{M}+1}$ so that
\begin{equation}\label{3*f6}
\text{if }t\not\equiv t' \text{ mod }2^{k_{M}+1}\text{ then }
(\pi(t),\pi(t'))=1.
\end{equation}
These periods can be different primes.
We suppose that if $\pi^{*}=\max_{t}\pi(t)$ then
\begin{equation}\label{*g7}
2^{n_{1}}\gg (\pi^{*})^{k_{M}}
\end{equation}
and
\begin{equation}\label{*f6}
\text{ if }x,x+\pi(t)\in [(t-1)2^{n_{1}}+1,t2^{n_{1}}+1)
\text{ then }\ppp_{k_{j}'}(x)=\ppp_{k_{j}'}(x+\pi(t)),
\end{equation}
this means that $\ppp_{k_{j}'}$ is ``periodic"
by $\pi(t)$ within $[(t-1)2^{n_{1}}+1,t2^{n_{1}}+1)$.
We will choose $\ppp_{k_{j}'}$ so that it is constant
on intervals of the form $[n,n+1)$ for all $n\in\N$,
it is periodic by $2^{n_{1}+k_{M}}$,
its range is $ \{ 2^{k_{j}'},0  \}$ and
\begin{equation}\label{2*f6}
\frac{1+\rrr}{2}\eee\cdot 2^{n_{1}}<
\int_{[(t-1)2^{n_{1}}+1,t2^{n_{1}}+1) }
\ppp_{k_{j}'}(x)dx <\eee 2^{n_{1}},
\end{equation}
provided $[(t-1)2^{n_{1}},t2^{n_{1}})\sse
\isss_{k_{j}}(\fff_{n_{1}})\sm \isss_{k_{j-1}}(\fff_{n_{1}})$.
Inequality \eqref{2*f6}
can be achieved if $2^{n_{1}}$ is sufficiently large.

We need to say something about the points where
\begin{equation}\label{*f7a}
\frac{\fff_{n_{1}}(x+l)\ppp_{k_{j}'}(x+2l)}{l}\geq 1.
\end{equation}
Since $\fff_{n_{1}}(x+l)=2^{n_{1}}$ if $x+l\in
\spt (\fff_{n_{1}})$ and
$\ppp_{k_{j}'}(x+2l)=2^{k_{j}'}$ if $x+2l\in
\spt (\ppp_{k_{j}'})$ we need to consider
$l\leq 2^{n_{1}+k_{j}'}.$

{\it A few words about our general plan.} We will be interested in
certain intervals $[t_{0}2^{n_{1}},(t_{0}+1)2^{n_{1}})\sse
\isss_{k_{j}}(\fff_{n_{1}})\sm \isss_{k_{j-1}}(\fff_{n_{1}})$, the
exact assumption about these intervals will be given in
\eqref{*f7c}. In these intervals
$[t_{0}2^{n_{1}},(t_{0}+1)2^{n_{1}})$ we consider sets
$\cax_{t_{0}}$ satisfying \eqref{3*f7h}. These sets are unions of
the subsets $\cax_{t_{0},t_{1}}$ see \eqref{*30*}. For the sets
$\cax_{t_{0},t_{1}}$ we have \eqref{*25x2}. For fixed $t_{0}$ but
different $t_{1}$'s the sets $\cax_{t_{0},t_{1}}$ are sufficiently
independent, so we have \eqref{*f7g}. Based on this we can obtain a
lower estimate of the measure of $\cax_{t_{0}}$ see  \eqref{2*f7g},
\eqref{3*f7g} and \eqref{*f7h}.

\subsection{SUBSTEP 1c: The auxiliary sets $X'(j,n_{1}+k_{j}')$,\\
$X(1/2,3/4,j)$, $X(1/2,3/4,j,e)$ and $X(1/2,3/4,j,o)$}\label{213}

We can suppose that $k_{j}\gg k_{j}'= k_{j-1}+10$ is so large that
for most points $x\in \isss_{k_{j}}(\fff_{n_{1}})\sm
\isss_{k_{j-1}}(\fff_{n_{1}})$ we have
$[x,x+2\cdot 2^{n_{1}+k_{j}'}+2^{n_{1}})\sse \isss_{k_{j}}(\fff_{n_{1}})$.
We denote the set of these $x$'s by $X'(j,n_{1}+k_{j}').$
Hence, if $k_{j}$ is sufficiently larger than
$k_{j}'$ then
\begin{equation}\label{2*f7a}
1\approx \frac{\olll(X'(j,n_{1}+k_{j}'))}
{\olll(\isss_{k_{j}}(\fff_{n_{1}})\sm
\isss_{k_{j-1}}(\fff_{n_{1}}))}>\frac{1}2.
\end{equation}
We remark that our construction implies that in \eqref{2*f7a} the choice of
$k_{j}$ does not depend on the choice of $n_{1}$, which means that one
can choose a $K_{j,0}$ such that for all $k_{j}\geq K_{j,0}$ we have
\eqref{2*f7a} for any $n_{1}$.
Denote by $X(1/2,3/4,j)$ the set of those $x$
for which there exists $t\in\Z$ such that
\begin{equation}\label{*g9}
x\in [(t+\frac12)2^{n_{1}+k_{j-1}},
(t+\frac34) 2^{n_{1}+k_{j-1}}).
\end{equation}
We split $X(1/2,3/4,j)$ into two subsets depending on
the parity of $t$.

If
\eqref{*g9} holds and
$t$ is even
then $x\in X(1/2,3/4,j,e)$,
while for odd $t$'s
$x\in X(1/2,3/4,j,o)$.

Suppose $x,y\in X(1/2,3/4,j)$, $x<y$, $x\in X(1/2,3/4,j,e)$,\\
$y\in X(1/2,3/4,j,o)$, and $l=\lf y\rf -\lf x\rf >0.$ Then, as the reader
can verify, $\lf x \rf\in X(1/2,3/4,j,e)$, $y'\defeq x+l\in X(1/2,3/4,j,o)$
and
$x+2l\in [2t'\cdot 2^{n_{1}+k_{j-1}},
(2t'+1)2^{n_{1}+k_{j-1}})$ with a $t'\in\Z$.
Hence,
\begin{align}\label{*f7b}
&\text{if }x\in X'(j,n_{1}+k_{j}')\cap X(1/2,3/4,j,e)
\text{ and }
\\&\nonumber
y\in X(1/2,3/4,j,o),\  0<l=\lf y\rf -\lf x\rf \leq 2^{n_{1}+k_{j}'}\text{ then}
\\&\nonumber
y'=x+l\in \isss_{k_{j-1}}(\fff_{n_{1}})\text{ and }
x+2l\in \isss_{k_{j}}(\fff_{n_{1}})\sm\isss_{k_{j-1}}(\fff_{n_{1}}).
\end{align}

\subsection{SUBSTEP 1d: Estimate of the measure of those points
where \eqref{*f7a} holds in one $2^{n_{1}}$ grid interval,
definition of the sets $\cax_{t_{0},t_{1}}$}

Suppose
\begin{equation}\label{*f7c}
[t_{0}2^{n_{1}},(t_{0}+1)2^{n_{1}})
\sse X'(j,n_{1}+k_{j}')\cap X(1/2,3/4,j,e).
\end{equation}
We want to obtain an estimate of the measure of those
$x$'s for which \eqref{*f7a} holds. By \eqref{*f7c}
there exists $t_{0}'\in\Z$ such that
\begin{equation}\label{4*f7ca}
[t_{0}2^{n_{1}},(t_{0}+1)2^{n_{1}})\sse
[(2t_{0}'+\frac{1}2)2^{n_{1}+k_{j-1}},
(2t_{0}'+\frac{3}4)2^{n_{1}+k_{j-1}}).
\end{equation}
Suppose
\begin{equation}\label{**f7c}
y\in X(1/2,3/4,j,o)\cap \supp(\fff_{n_{1}})
\text{ with }0<l=\lf y\rf -\lf x\rf \leq 2^{n_{1}+k_{j}'}.
\end{equation}
Then $\fff_{n_{1}}(y)=2^{n_{1}}$ and by \eqref{*f5}
\begin{equation}\label{4*f7c}
\text{there exists }t_{1}\in\Z\text{ such that }
y\in[t_{1}2^{n_{1}},t_{1}2^{n_{1}}+1).
\end{equation}
Moreover, $l=t_{1}2^{n_{1}}-\lf x \rf$ and $y'=x+l\in [t_{1}2^{n_{1}},t_{1}2^{n_{1}}+1)$
as well which implies $\fff_{n_{1}}(x+l)=\fff_{n_{1}}(y')=2^{n_{1}}$
and $|y-y'|<1.$

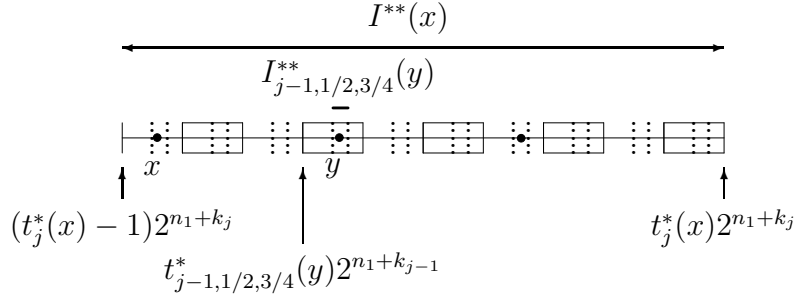
\placedrawing[h]{fgt3.lp}{Notation related to $I_{j}^{**}(x)$}{figure3}

Denote by
\begin{align}\label{3*f7c}
\Ics j (x)=&[(\tcs j (x)-1)2^{n_{1}+k_{j}},
\tcs j (x)2^{n_{1}+k_{j}})\\&\nonumber
\text{the component of
 }\isss_{k_{j}}(\fff_{n_{1}})\text{ containing }
x.
\end{align}
By the definition of $X'(j,n_{1}+k_{j}')$ the points
$y'=x+l,$ and $x+2l$ belong to $\Ics j (x)$, that is,
$x,$ $y'=x+l$, and $x+2l$ belong to the same component
of $\isss_{k_{j}}(\fff_{n_{1}})$.
Moreover, by \eqref{*f7b}, $y'=x+l\in\isss_{k_{j-1}}
(\fff_{n_{1}})$ and $x+2l\in\isss_{k_{j}}(\fff_{n_{1}})
\sm \isss_{k_{j-1}}(\fff_{n_{1}}).$
We introduce the notation
$$\Ics{j-1,1/2,3/4}(y)=
[(\tcs{j-1,1/2,3/4}(y)+\frac{1}2)2^{n_{1}+k_{j-1}},
(\tcs{j-1,1/2,3/4}(y)+\frac{3}4)2^{n_{1}+k_{j-1}})$$
for the interval containing $y\in X(1/2,3/4,j)$.
Then $y'\in\Ics{j-1,1/2,3/4}(y)$ holds as well.
On Figure \ref{figure3}, $\Ics{j-1,1/2,3/4}(y)$ is the tiny interval containing
$y$, its length is marked by a short line segment above $y$.
Keep in mind that we supposed that $j\geq 2$.
By our construction
if
\begin{align}\label{3*f7d}
I^{*} =&[(t^{*}+\frac{1}2)2^{n_{1}+k_{j-1}},(t^{*}+\frac{3}4)
2^{n_{1}+k_{j-1}})\text{ is a component of }\\\nonumber
&\isss_{k_{j-1}}(\fff_{n_{1}})\cap
X(1/2,3/4,j)\text{ then }
\\ \label{*f7d}
&\lll(\isss_{0}(\fff_{n_{1}})\cap I^{*})=
\bigg(\frac{1}2\bigg)^{j-1}\lll(I^{*})=
\bigg(\frac{1}2\bigg)^{j-1}2^{n_{1}+k_{j-1}-2}.
\end{align}
This implies that there exist exactly
\begin{align}\label{**f7d}
&\bigg(\frac{1}2\bigg)^{j-1}2^{k_{j-1}-2}\text{ many }t_{1}
\text{'s such that }\\
\nonumber &
[t_{1}2^{n_{1}},t_{1}2^{n_{1}}+1)\sse\supp(\fff_{n_{1}})
\cap I^{*}.
\end{align}
We denote the set of these $t_{1}$'s by $\cat_{1}(I^{*})$.

We still assume that $t_{0}$ satisfies \eqref{*f7c}.
By \eqref{*isptdef}
in the interval $[t_{0}2^{n_{1}}+2^{n_{1}+k_{j}'-1},
t_{0}2^{n_{1}}+2^{n_{1}+k_{j}'})$
there are $2^{k_{j}'-k_{j-1}-2}$ many $I^{*}$'s satisfying
\eqref{3*f7d}. Denote the set of the corresponding $t^{*}$'s
by $\cat^{*}(t_{0})$. Finally, denote by $\cat_{1}(t_{0})$
the set of those $t_{1}$ which belong to a $\cat_{1}(I^{*})$
with $I^{*}$ of the form in \eqref{3*f7d} and $t^{*}\in
\cat^{*}(t_{0})$.
Then by \eqref{**f7d}
\begin{equation}\label{*f7e}
\# \cat_{1}(t_{0})=\left (\frac{1}2\right)^{j-1}2^{k_{j-1}-2}
2^{k_{j}'-k_{j-1}-2}
=
\left (\frac{1}2\right)^{j}
2^{k_{j}'-3}.
\end{equation}

Next suppose $t_{1}\in\cat_{1}(t_{0})$ is fixed.
If $y\in [t_{1}2^{n_{1}},t_{1}2^{n_{1}}+1)$
then by \eqref{**f7d} $\fff_{n_{1}}(y)=2^{n_{1}}$.
For $x\in [t_{0}2^{n_{1}},(t_{0}+1)2^{n_{1}})$
set $l_{x,t_{1}}=t_{1}2^{n_{1}}-\lf x \rf=\lf y \rf-\lf x \rf$
and $y'=x+l_{x,t_{1}}.$
Then $\fff_{n_{1}}(x+l_{x,t_{1}})=\fff_{n_{1}}(y')=\fff_{n_{1}}(y)=2^{n_{1}}$.
From $t_{1}\in \cat_{1}(t_{0})$ it follows that
$[t_{1}2^{n_{1}}, t_{1}2^{n_{1}}+1)\sse [t_{0}2^{n_{1}}+2^{n_{1}+k_{j}'-1},
t_{0}2^{n_{1}}+2^{n_{1}+k_{j}'})$. This and $x\in [t_{0}2^{n_{1}},(t_{0}+1)2^{n_{1}})$
implies
\begin{equation}\label{*38bb*}
0<l_{x,t_{1}}=t_{1}2^{n_{1}}-\lf x \rf\leq (t_{1}-t_{0})2^{n_{1}}
<2^{n_{1}+k_{j}'}<2^{n_{1}+k_{M}}.
\end{equation}
We have
$\fff_{n_{1}}(x+l_{x,t_{1}})=2^{n_{1}}$
and
\begin{equation}
\label{*25x}
\frac{\fff_{n_{1}}(x+l_{x,t_{1}})\ppp_{k_{j}'}(x+2l_{x,t_{1}})}
{l_{x,t_{1}}}\geq 1
\end{equation}
 if $x+2l_{x,t_{1}}\in\supp(\ppp_{k_{j}'})$.
For $t_{1}\in \cat_{1}(t_{0})$
denote by $\cax_{t_{0},t_{1}}$
the set of those $x\in [t_{0}2^{n_{1}},(t_{0}+1)2^{n_{1}})$
for which $x+2l_{x,t_{1}}\in\supp(\ppp_{k_{j}'}).$
This means that
\begin{equation}\label{*25x2}
\text{for }x\in \cax_{t_{0},t_{1}}\sse [t_{0}2^{n_{1}},(t_{0}+1)2^{n_{1}})
\text{ there exists }l_{x,t_{1}}\text{ such that}
\end{equation}
$$\text{\eqref{*25x} holds and }y'=x+l_{x,t_{1}}\in [t_{1}2^{n_{1}},
t_{1}2^{n_{1}}+1).$$

\subsection{SUBSTEP 1e: ``Periodicity and independence"
of the sets $\cax_{t_{0},t_{1}}$}

\placedrawing[h]{fgt6.lp}{Sets $\cax_{t_{0,t_{1}}}$ and $\cax_{t_{0},t_{1}'}$}
{figure6}

Observe that if $x\in [t_{0}2^{n_{1}},(t_{0}+1)2^{n_{1}})$
then
$\lf x\rf \in [t_{0}2^{n_{1}},(t_{0}+1)2^{n_{1}}-1]$
and
$x+2l_{x,t_{1}}\in [(2t_{1}-t_{0}-1)2^{n_{1}}+1,
(2t_{1}-t_{0})2^{n_{1}}+1)$, moreover by \eqref{*f6}
$\spt (\ppp_{k_{j}'})$ is ``periodic" by $\pi(2t_{1}-t_{0})$,
this implies that
\begin{align}\label{*f7f}
&\text{if }x,x+\pi(2t_{1}-t_{0})\in[t_{0}2^{n_{1}},(t_{0}+1)
2^{n_{1}})\text{ then }\\&
\nonumber x\in
\cax_{t_{0},t_{1}}
\text{ iff }x+\pi(2t_{1}-t_{0})\in
\cax_{t_{0},t_{1}}.
\end{align}

From \eqref{2*f6} it follows that
\begin{equation}\label{2*f7f}
\frac{1+\rrr}2
\eee2^{-k_{j}'}\cdot 2^{n_{1}}<
\lll(\cax_{t_{0},t_{1}})<
\eee2^{-k_{j}'}\cdot 2^{n_{1}}.
\end{equation}
In \eqref{*f7g} below it will be useful to keep in mind
that the ``density" of $\cax_{t_{0},t_{1}'}$
in $[t_{0}2^{n_{1}},(t_{0}+1)2^{n_{1}})$
is $\lll(\cax_{t_{0},t_{1}'})/2^{n_{1}}$ and
\eqref{2*f7f} holds for $t_{1}'$ as well.
By \eqref{3*f6}, \eqref{*38bb*}
and \eqref{*f7f} for $t_{1},t_{1}'\in\cat_{1}(t_{0})$,
$t_{1}\not=t_{1}'$ the ``periods" of
$\cax_{t_{0},t_{1}}$ and $\cax_{t_{0},t_{1}'}$ are
relatively prime, see Figure \ref{figure6}.
By Lemma \ref{*indep} if $n_{1}$ is sufficiently large
 these sets are independent in the
sense that
\begin{equation}\label{*f7g}
\lll(\cax_{t_{0},t_{1}}\cap
\cax_{t_{0},t_{1}'})<\lll(\cax_{t_{0},t_{1}})2
\frac{\lll(\cax_{t_{0},t_{1}'})}{2^{n_{1}}}<
2\eee2^{-k_{j}'}\lll(\cax_{t_{0},t_{1}}).
\end{equation}
Hence, using \eqref{*f7e} and \eqref{*f7g} we infer
\begin{align}\label{2*f7g}
&\lll(\cax_{t_{0},t_{1}}\sm\bigcup_{t_{1}'\not=t_{1},\
 t_{1}'\in\cat_{1}(t_{0})}\cax_{t_{0},t_{1}'})>\\
\nonumber &
(1-2\eee\cdot 2^{-k_{j}'}\left (\frac{1}2\right)^{j}
2^{k_{j}'-3})\lll(\cax_{t_{0},t_{1}})>
\frac{1}2\lll(\cax_{t_{0},t_{1}}).
\end{align}
This, \eqref{*f7e}, and \eqref{2*f7f} imply that
\begin{align}\label{3*f7g}
&\lll(\bigcup_{t_{1}\in\cat_{1}(t_{0})}\cax_{t_{0},t_{1}})>
\frac{1}2\sum_{t_{1}\in\cat(t_{0})}\lll(\cax_{t_{0},t_{1}})>\\
&\nonumber
\frac{1}4\eee2^{-k_{j}'}\cdot 2^{n_{1}}\left (\frac{1}2\right)^{j}
\cdot 2^{k_{j}'-3}=\frac{\eee}{32}2^{n_{1}}\left (\frac{1}2\right)^{j}.
\end{align}
This means that in each interval
$I(t_{0})=[t_{0}2^{n_{1}},(t_{0}+1)2^{n_{1}})$ satisfying
\eqref{*f7c} we could find a set
\begin{equation}\label{*30*}
\cax_{t_{0}}=\bigcup_{t_{1}\in\cat_{1}(t_{0})}\cax_{t_{0},t_{1}},
\end{equation}
such that
\begin{equation}\label{*f7h}
\lll(\cax_{t_{0}})>\frac{\eee}{32} \left (\frac{1}2\right)^{j}
\lll(I(t_{0}))
\end{equation}
and
\begin{equation}\label{3*f7h}
\text{for }x\in\cax_{t_{0}},\  \exists l\leq 2^{n_{1}+k_{M}}
 \text{ such that }
\frac{\fff_{n_{1}}(x+l)\ppp_{k_{j}'}(x+2l)}{l}\geq 1.
\end{equation}

\subsection{SUBSTEP 1f: Conclusion of the $s=1$ case}

Denote by $\cai(j)$  the union of all intervals in
$\isss_{k_{j}}(\fff_{n_{1}})\sm \isss_{k_{j-1}}(\fff_{n_{1}})$
which satisfy \eqref{*f7c}.
Using
\eqref{*ISP0*}, \eqref{*ISP*}, \eqref{2*f7a} and \eqref{*g9} we obtain
\begin{equation}\label{2*f7h}
\olll(\cai(j))>\frac{1}{16}\olll(\isss_{k_{j}}(\fff_{n_{1}}))=
\frac{1}{16}2^{-M+j}.
\end{equation}
Denote by $\ocax(j)$ the set of those $x\in\isss_{k_{j}}(\fff_{n_{1}})
\sm \isss_{k_{j-1}}(\fff_{n_{1}})$ for which there exists
$l$ such that \eqref{3*f7h} holds. Then
\eqref{*f7h} and \eqref{2*f7h} imply that
\begin{equation}\label{*f7i}
\olll(\ocax(j))>\frac{\eee}{512}\left (\frac{1}2\right)^{j}
\cdot 2^{-M+j}=\frac{\eee}{512}2^{-M}.
\end{equation}
The sets $\ocax(j)\sse \isss_{k_{j}}(\fff_{n_{1}})\sm\isss_{k_{j-1}}
(\fff_{n_{1}})$ and
$\ocax(j')\sse \isss_{k_{j'}}(\fff_{n_{1}})\sm\isss_{k_{j'-1}}
(\fff_{n_{1}})$ are disjoint when $j\not=j'$ and hence using
\eqref{*f7i} for $j=2,...,M$ one can see that \eqref{*f1}
holds with $\kkks=1.$ We set $\ooo=k_{M}'.$
From \eqref{*AVINTPHIN1*} it follows \eqref{**f1}
when $s=1$.
From  \eqref{*psiksup1*}, \eqref{2*f6b}, \eqref{*psiksup}
and \eqref{2*f6} it follows
\eqref{***f1}.

\subsection{STEP 2: The general step of the induction\\
SUBSTEP 2a: usage of the functions from step $s-1$
of the induction}\label{subsubstep2a}

Next we turn to the general step of our induction. Suppose that for
any possible choice of $\widetilde{\aaa},$ $\widetilde{\capp}$,
$\widetilde{\rrr},$ $0<\widetilde{\eee}<\eee_{0}$, and
$\widetilde{M}$ one can find
$\widetilde{\aaa}-\widetilde{\capp}-(\kkks-1)$-families. In case in
\eqref{*f1}, $\olll \{ .  \}>\eee$ holds for an $\aaa-
\capp-(\kkks-1)$-family then one can define an
$\aaa-\capp-\kkks$-family by choosing an almost arbitrary
$\fff_{n_{\kkks-1}+1}= \fff_{n_{\kkks}}$ so that \eqref{**f1} holds.
So we can assume that we work with $\aaa-\capp-(\kkks-1)$-families
for which in \eqref{*f1}, $\olll \{ .  \}\leq\eee$ holds. We define
$\fff_{n_{\kkks}}$ analogously to $\fff_{n_{1}}$. The interval
supports are defined by
\begin{align}\label{*8ns}
\isss_{k}(\fff_{n_{s}})=  \bigcup  &\{ [(t-1)2^{n_{s}+k},t2^{n_{s}+k})
:\\
\nonumber
&t\in \Z,\   \spt(\fff_{n_{s}})\cap  [(t-1)2^{n_{s}+k},t2^{n_{s}+k})
\not=\ess \} .
\end{align}

We set $$\isss_{k_{M}}(\fff_{n_{s}} )=\R.$$
\begin{align}\nonumber
&\text{For }j\in  \{ 1,...,M  \} \text{ an interval }
[(t-1)2^{n_{s}+k_{j-1}},t 2^{n_{s}+k_{j-1}})\sse
\isss_{k_{j}}(\fff_{n_{s}} )
 \\
&\label{*isptdefns}\text{belongs to }
\isss_{k_{j-1}}(\fff_{n_{s}} ) \text{ if and only if }
t \text{ is even.}\end{align}
This implies
\begin{equation}\label{*9ns*}
\olll(\isss_{k_{j}}(\fff_{n_{s}}))=\left (\frac{1}{2} \right )^{M-j},
\end{equation}
and
\begin{equation}\label{*10ns*}
\olll(\isss_{k_{j-1}}(\fff_{n_{s}}))=\olll(\isss_{k_{j}}(\fff_{n_{s}})\sm
\isss_{k_{j-1}}(\fff_{n_{s}}))=\frac{1}{2}\olll(\isss_{k_{j}}(\fff_{n_{s}})).
\end{equation}

 It will be useful to keep in mind
for further reference that $t-1$, which corresponds
to the left endpoint of the support intervals is odd.

We will define $\fff_{n_{s}} $ so that it will be periodic by
$2^{n_{s}+k_{M}}\in \capp$.

If $[(t-1)2^{n_{s}},t 2^{n_{s}})\sse\isss_{0}(\fff_{n_{s}})=\isss_{k_{0}}(\fff_{n_{s}})$
then
\begin{align}\label{*11ns}
&\fff_{n_{s}}(x)=2^{n_{s}}\text{ if }x\in [(t-1)2^{n_{s}},(t-1)2^{n_{s}}
+1),\text{ and }\\
&\fff_{n_{s}}(x)=0\text{ if }x\in [(t-1)2^{n_{s}}+1,t\cdot 2^{n_{s}}
).
\end{align}

From \eqref{*9ns*} and \eqref{*10ns*} used with $j=1$ it follows that
\begin{equation}\label{*13ns*}
\ointt\fff_{n_{s}}=2^{-M}.
\end{equation}

We need more assumptions about $k_{0}\ll...\ll k_{M}$. Set
$\aaa_{1}=\aaa$. Choose
 an
$\aaa_{1}-\capp-(\kkks-1)$-family with $\ppp$-interval
$[\aaa_{1},\ooo_{1}]$ periodic by $p_{1,\kkks-1}\in \capp$. Recall
that $\capp $ consists of powers of $2$ and hence $p_{1,\kkks-1}$ is
also a power of $2$. As we remarked earlier we can suppose that for
this family in \eqref{*f1}, $\olll \{ .  \}\leq\eee$ holds. We
denote the functions corresponding to this family by
$\fff\so1,...,\fff\so{n_{1,\kkks-1}}$,
$\ppp\so{\aaa_{1}},...,\ppp\so{\ooo_{1}}$. For other parameters
belonging to this family we adopt a similar subscript notation. We
choose $k_{1}$ so that $k_{1}\gg p\so{\kkks-1}\geq\ooo_{1}.$ If
$x\in\isss_{k_{1}}(\fff_{n_{\kkks}})$ then we set
$\fff_{i}(x)=\fff\so i(x)$ for $i=1,...,n\so{\kkks-1}$ and
$\ppp_{k}(x)=\ppp\so k(x)$ for $k=\aaa_{1},...,\ooo_{1}$, moreover,
for other $k$'s we set $\ppp_{k}(x)=0$ and for
$i=n_{1,s-1}+1,...,n_{s}-1$ we set $\fff_{i}(x)=0$. If $k_{1}$ is
sufficiently larger than $p\so{\kkks-1}$ then the length of the
 components
of $\isss_{k_{1}}(\fff_{n_{\kkks}})$
will be a multiple of $p_{1,s-1}$
 and hence \eqref{**f1}
and \eqref{***f1} will stay valid ``relative to"
$\isss_{k_{1}}(\fff_{n_{\kkks}})$.
By this we mean the following:
\begin{equation}\label{*3rel}
(s-1)\cdot 2^{-M-1}\olll(\isss_{k_{1}}(\fff_{n_{s}}))<\sum_{i=1}^{n_{s}-1}\ointt\fff_{i}
\restr{\isss_{k_{1}}(\fff_{n_{s}})}<
\end{equation}
$$(\kkks-1)\cdot 2^{-M+1}\olll(\isss_{k_{1}}(\fff_{n_{s}})),$$
and
\begin{equation}\label{*4rel}
\rrr\eee\olll(\isss_{k_{1}}(\fff_{n_{s}}))<\sum_{k=\aaa_{1}}^{\ooo_{1}}\ointt
\ppp_{k}\restr{\isss_{k_{1}}(\fff_{n_{s}})}
=
\end{equation}
$$\sum_{k=\aaa}^{\ooo}\ointt
\ppp_{k}\restr{\isss_{k_{1}}(\fff_{n_{s}})}<\eee
\olll(\isss_{k_{1}}(\fff_{n_{s}})).$$

Suppose $j\geq 2$ and $k_{j-1}$ is defined. Let
$\aaa_{j}=k_{j-1}+10$ and choose an
$\aaa_{j}-\capp-(\kkks-1)$-family with $\ppp$-interval
$[\aaa_{j},\ooo_{j}].$ Again, as we remarked earlier we can suppose
that for this family $\olll \{ .  \}\leq \eee$ holds in \eqref{*f1}.
Denote the corresponding functions by
$\fff\sj1,...,\fff\sj{n_{j,\kkks-1}}$,
$\ppp\sj{\aaa_{j}},...,\ppp\sj{\ooo_{j}}$ periodic by
$p\sj{\kkks-1}$. We will choose $k_{j}\gg p\sj{\kkks-1}\geq
\ooo_{j}$.

We repeat the above steps for $j=2,...,M$ and obtain
\begin{equation}\label{*40b*}
k_{M}\gg p_{M,s-1}\geq \ooo_{M}\geq\aaa_{M}\gg
k_{M-1}\gg ...
\end{equation}
 $$ \gg k_{j}\gg p_{j,s-1}\geq \ooo_{j}
\geq \aaa_{j}\gg k_{j-1}\gg...\gg k_{1}\gg p_{1,s-1}\geq\ooo_{1}.$$
The $\ppp$-interval of our $\aaa-\capp-\kkks$-family will be defined
by $\aaa=\aaa_{1}$ and $\ooo=\ooo_{M}$.

In the end we choose $n_{\kkks}$ so large that
\begin{equation}\label{*g20}
n_{\kkks}>\max\{ n\sj{\kkks-1}: j=1,...,M  \}.
\end{equation}
Similarly to the $s=1$ case our parameters are chosen in the order
$\aaa=\aaa_{1}$, $\ooo_{1}$, $p_{1,s-1}$, $k_{1}$, ..., $k_{j-1}$,
$\aaa_{j}$, $\ooo_{j}$, $p_{j,s-1}$, $k_{j}$, ..., $k_{M}$, and then
we choose some new parameters, the $\pi(t,j)$'s and finally we fix a
large $n_{s}$. At the $j$'th step we have a $k_{j-1}$ which
determines $\aaa_{j}= k_{j-1}+10$. The $\aaa_{j}-\capp-(s-1)$ family
provides $\ooo_{j}$ and $p_{j,s-1}$. If we have the value of
$p_{j,s-1}$ then we choose $k_{j}$ so that \eqref{*f11} holds. We
emphasize again that, similarly to \eqref{2*f7a}, \eqref{*f11} does
not depend on the choice of $n_{s}$.

Since $\capp$ consists of powers of $2$
all $p\sj{\kkks-1}$ are powers of $2$.
To each $j=1,...,M$ and each interval of the form
$[(t-1)2^{n_{\kkks}},t2^{n_{\kkks}})$, ($t\in\Z$) we assign
a set $\capp(t,j)$ consisting of infinitely
many odd numbers such that if $\pi(t,j)\in\capp(t,j)$
then
\begin{equation}\label{*f9}
\pi(t,j)\gg 2^{k_{M}}p_{M,\kkks-1}
\end{equation}
and
\begin{align}\label{2*f9}
&\text{if }t\not\equiv t'\text{ mod }2^{k_{M}+1},\ \pi(t,j)\in\capp(t,j),
\pi(t',j)\in\capp(t',j)\\
&\nonumber \text{then }(\pi(t,j),\pi(t',j))=1.
\end{align}
On the other hand,
\begin{equation}\label{3*f9}
\text{if }t\equiv t'\text{ mod }2^{k_{M}+1}\text{ then }
\capp(t,j)=\capp(t',j),
\end{equation}
moreover
\begin{equation}\label{*g21}
(\pi(t,j),\pi(t',j'))=1\text{ if }j\not=j',\
\pi(t,j)\in\capp(t,j),\
\pi(t',j')\in\capp(t',j')
.
\end{equation}
We have already defined for $x\in\isss_{k_{1}}(\fff_{n_{\kkks}})$
the values $\fff_{i}(x)$ for $i=1,...,n_{\kkks}$ and
$\ppp_{k}(x)$ for all $k.$

We also want to define these functions when
$x\not\in\isss_{k_{1}}(\fff_{n_{\kkks}})$, that is,
there exists a $j\geq 2$ such that
$x\in\isss_{k_{j}}(\fff_{n_{\kkks}})\sm \isss_{k_{j-1}}(\fff_{n_{\kkks}})$.

Suppose $j\geq 2.$ For any $t\in\Z$ from $\capp(t,j)$ choose and fix a number
$\pi(t,j)$ such that by Lemma \ref{lemf3} applied to the
$\aaa_{j}-\capp-(\kkks-1)$-family with $\ppp$-interval
$[\aaa_{j},\ooo_{j}]$ we can select an
$\aaa_{j}-\capp(t,j)-(\kkks-1)$-family with $\ppp$-interval
$[\aaa_{j},\ooo_{j}]$ periodic by $\pi(t,j)$. Denote the
corresponding functions by $\fff\stj 1,...,\fff\stj{n_{j,\kkks-1}}$,
$\ppp\stj{\aaa_{j}},...,\ppp\stj{\ooo_{j}}.$ By \eqref{3*f9} we can
suppose that
\begin{equation}\label{3*f9b}
\text{if }t\equiv t'\text{ mod }2^{k_{M}+1}\text{ then }
\pi(t,j)=\pi(t',j),
\end{equation}
$$\fff\stj i=\fff_{t',j,i}\text{ for }i=1,...,n_{j,s-1},$$
$$\text{and } \ppp\stj k=\ppp_{t',j,k}\text{ for }k=\aaa_{j},...,\ooo_{j}.$$
We also have \eqref{*f9}. By \eqref{2*f9} and \eqref{*g21}
\begin{align}\label{2*f9b}
&\text{if }t\not\equiv t'\text{ mod }2^{k_{M}+1}\text{ then }
(\pi(t,j),\pi(t',j))=1, \text{ moreover}\\
&\nonumber
(\pi(t,j),\pi(t',j'))=1 \text{ if }j\not=j',\text{ for any choice of }
t\text{ and }t'.
\end{align}
Based on \eqref{*f16} we can suppose that
the $\pi(t,j)$'s are chosen so large that
\begin{equation}\label{2*f16}
\ointt\ppp\stj k\approx\ointt \ppp_{j,k},\text{ for all }t,j,k.
\end{equation}
By using $\rrr\eee<\sum_{k=\aaa_{j}}^{\ooo_{j}}\ointt\ppp_{j,k}<\eee$
and \eqref{2*f16} we can obtain for large $\pi(t,j)$'s
\begin{equation}\label{3*f16}
\rrr\eee<\sum_{k=\aaa_{j}}^{\ooo_{j}}\min_{t}\ointt\ppp_{t,j,k}
\leq
\sum_{k=\aaa_{j}}^{\ooo_{j}}\max_{t}\ointt\ppp_{t,j,k}<\eee.
\end{equation}
Similarly,
based on \eqref{*5*5} we can suppose that
the $\pi(t,j)$'s are chosen so large that
\begin{equation}\label{2*f16b}
\ointt\fff\stj i\approx\ointt \fff_{j,i},\text{ for all }t,j,i.
\end{equation}
By using $(s-1)2^{-M-1}<\sum_{i=1}^{n_{s-1}}\ointt\fff_{j,i}<(s-1)2^{-M+1}$
and \eqref{2*f16b} we can obtain  for large $\pi(t,j)$'s
\begin{equation}\label{*57PPP}
(s-1)2^{-M-1}<\sum_{i=1}^{n_{j,s-1}}\min_{t}\ointt\fff_{t,j,i}
\leq
\sum_{i=1}^{n_{j,s-1}}
\max_{t}\ointt\fff_{t,j,i}<(s-1)2^{-M+1}.
\end{equation}

We suppose that
\begin{equation}\label{*48b*}
\text{if }\pi^{*}=\max_{t,j}\pi(t,j)
\text{ then }2^{n_{\kkks}}\gg (\pi^{*})^{k_{M}}
\end{equation}
(this means that $n_{\kkks}$ should be sufficiently large).

\subsection{SUBSTEP 2b: Definitions of $\fff_{i}$,
$i=1,...,n_{j,s-1}$ and $\ppp_{k}$, $k=\aaa_{j},...,\ooo_{j}$}

Suppose
\begin{equation}\label{*73oo*}
x\in [(t-1)2^{n_{\kkks}}+1,t2^{n_{\kkks}}+1)\cap
(\isss_{k_{j}}(\fff_{n_{\kkks}}) \sm\isss_{k_{j-1}}
(\fff_{n_{\kkks}})).
\end{equation}

Set
\begin{align}\label{*44b}
&\fff_{i}(x)=\fff\stj i(x) \text{ for }i=1,...,n\sj{\kkks-1},
\text{ and }\\
\nonumber
&\ppp_{k}(x)=\ppp\stj k(x) \text{ for }k=\aaa_{j},...,\ooo_{j}.
\end{align}
For other $k$'s set $\ppp_{k}(x)=0.$
\begin{equation}\label{*73xx*}
\text{We also put }\fff_{i}(x)=0\text{ for }n\sj{\kkks-1}<i<n_{\kkks}.
\end{equation}
See Figure \ref{figure4}. It is also clear that by the definition
of interval supports we have
\begin{equation}\label{*73zz*}
\fff_{n_{s}}(x)=0
\text{ for any }x\in \isss_{k_{j}}(\fff_{n_{\kkks}})\sm
\isss_{k_{j-1}}(\fff_{n_{\kkks}}).
\end{equation}

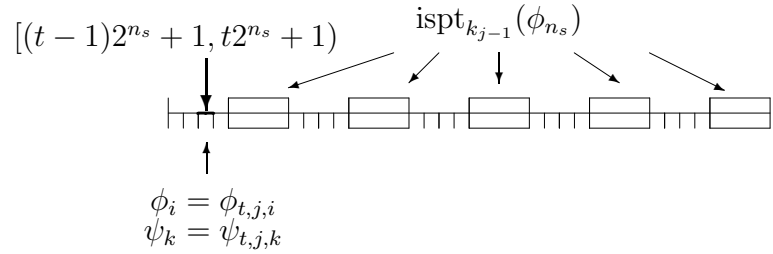
\placedrawing[h]{fgt4.lp}{Definition of $\fff_{i}$ and $\ppp_{k}$
in one component of $\isss_{k_{j}}(\fff_{n_{s}})$}{figure4}

By the periodicity assumptions about our corresponding $\kkks-1$
families we have
\begin{align}\label{*f9c}
&\text{if }x,x+\pi(t,j)\in [(t-1)2^{n_{\kkks }}+1,
t2^{n_{\kkks}}+1)\cap\\
&\nonumber  (\isss_{k_{j}}(\fff_{n_{\kkks}})\sm
\isss_{k_{j-1}}(\fff_{n_{\kkks}}))
\text{ then }\ppp_{k}(x)=\ppp_{k}(x+\pi(t,j)).
\end{align}

If $n_{s}$ is sufficiently large then by our induction step about
$\aaa_{j}-\capp(t,j)-(\kkks-1)$-families we have
\begin{align}\label{*1043b}
\olll \bigg  \{ x: \max_{k\in[\aaa_{1},\ooo_{M}]}
\max_{l\leq \pi^{*}} &
\frac{\sum_{i=1}^{n_{\kkks-1}}\fff_{i}(x+l)\ppp_{k}(x+2l)}{l}\geq 1
 \bigg \} >\\ \nonumber
&
\frac{(\kkks-1)(M-1)\eee2^{-M}}{2048}.
\end{align}

By choosing $n_{s}$ sufficiently large
and keeping in mind \eqref{*73xx*}, \eqref{*73zz*}
we have a version
of \eqref{2*f16b} and \eqref{*57PPP} with respect to
$\isss_{k_{j}}(\fff_{n_{s}})\sm \isss_{k_{j-1}}(\fff_{n_{s}})$
\begin{equation}\label{**3rel}
(s-1)\cdot 2^{-M-1}\olll(\isss_{k_{j}}(\fff_{n_{s}})\sm \isss_{k_{j-1}}(\fff_{n_{s}}))<
\sum_{i=1}^{n_{s}-1}\ointt\fff_{i}
\restr{\isss_{k_{j}}(\fff_{n_{s}})\sm \isss_{k_{j-1}}(\fff_{n_{s}})}<
\end{equation}
$$(\kkks-1)\cdot 2^{-M+1}\olll(\isss_{k_{j}}(\fff_{n_{s}})\sm \isss_{k_{j-1}}(\fff_{n_{s}})),$$
and
\begin{equation}\label{**4rel}
\rrr\eee\olll(\isss_{k_{j}}(\fff_{n_{s}})\sm \isss_{k_{j-1}}(\fff_{n_{s}}))<\sum_{k=\aaa_{1}}^{\ooo_{1}}\ointt
\ppp_{k}\restr{\isss_{k_{j}}(\fff_{n_{s}})\sm \isss_{k_{j-1}}(\fff_{n_{s}})}
=
\end{equation}
$$\sum_{k=\aaa}^{\ooo}\ointt
\ppp_{k}\restr{\isss_{k_{j}}(\fff_{n_{s}})\sm \isss_{k_{j-1}}(\fff_{n_{s}})}<\eee
\olll(\isss_{k_{j}}(\fff_{n_{s}})\sm \isss_{k_{j-1}}(\fff_{n_{s}})).$$

\subsection{SUBSTEP 2c: Definitions related to $\fff_{n_{s}}$,
the sets $X'(j,n_{s}+\ooo_{j})$, $X(1/2,3/4,j)$, $X(1/2,3/4,j,e)$
and\\ $X(1/2,3/4,j,o)$}

We need to say something about the set of those $x$'s where
\begin{equation}\label{*f10}
\max_{k'\in[\aaa_{j},\ooo_{j}]} \frac{\fff_{n_{\kkks}}(x+l)\ppp_{k'}(x+2l)}
{l}\geq 1.
\end{equation}
We split the above set of $x$'s into subsets of those $x$'s
where
\begin{equation}\label{2*f10}
\frac{\fff_{n_{\kkks}}(x+l)\ppp_{k'}(x+2l)}
{l}\geq 1, \text{ with }\ooo_{j}\geq k'\geq \aaa_{j}> k_{j-1}.
\end{equation}
Since $\fff_{\nk}(x+l)=2^{\nk}$ if
$x+l\in \supp(\fff_{\nk})$ and $\ppp_{k'}(x+2l)=2^{k'}$
if $x+2l\in\supp(\ppp_{k'})$ we need to consider $l\leq 2^{\nk+k'}$.

{\it More words about our general plan.} We will be interested in
certain intervals $[t_{0}2^{n_{s}},(t_{0}+1)2^{n_{s}})\sse
\isss_{k_{j}}(\fff_{n_{s}})\sm \isss_{k_{j-1}}(\fff_{n_{s}})$, the
exact assumption about these intervals will be given in
\eqref{*f12}. In these intervals
$[t_{0}2^{n_{s}},(t_{0}+1)2^{n_{s}})$ we consider the sets
$\cax_{t_{0}}'$ satisfying \eqref{3*f7hf23}. By \eqref{*111*} these
sets will be the unions of the auxiliary sets $\cax'_{t_{0},k'}$
defined in \eqref{*83plus*} via the sets $\cax_{t_{0},t_{1},k'}$ and
$\cax_{t_{0},s-1}$. For the sets $\cax_{t_{0},t_{1},k'}$ with
$\aaa_{j}\leq k'\leq \ooo_{j}$ we have \eqref{*62*2}. For fixed
$t_{0}$ and $t_{1}$ by \eqref{*62*lxt} there can be at most one $k'$
for which an $\cax_{t_{0},t_{1},k'}$ is defined. For fixed $t_{0}$
but different $t_{1}$'s and $k'$'s the sets $\cax_{t_{0},t_{1},k'}$
are sufficiently independent, so we have \eqref{*f7gf17} and
\eqref{2*f7gf17}. It is a new feature that we also need to consider
the sets $\cax_{t_{0},s-1}$ defined in \eqref{**19*}.
 These
sets take care of points coming from the
$(s-1)$-families of earlier steps of the induction.
They will also be sufficiently independent of
the sets $\cax_{t_{0},t_{1},k'}$, see \eqref{*675} and then
(\ref{*f21}-\ref{2*f22}).
Hence, in \eqref{*f7hf23} we obtain our lower
estimate of the measure of the sets $\cax'_{t_{0}}$.
By \eqref{*87plus*} these sets are disjoint from $\cax_{t_{0},s-1}.$

Now we return to the details of the proof of Theorem \ref{f2}. One
can assume similarly to the one family case the following: Suppose
that $k_{j}\gg\ooo_{j}\geq k'$ is so large that for most points of
$x\in \isss_{k_{j}}(\fff_{n_{\kkks}})\sm
\isss_{k_{j-1}}(\fff_{n_{\kkks}}) $ we have $[x-2^{n_{s}},x+2\cdot
2^{\nk+\ooo_{j}}+2^{\nk})\sse \isss_{k_{j}}(\fff_{\nk})$. We denote
the set of these $x$'s by $X'(j,\nk+\ooo_{j})$. If
$k_{j}\gg\ooo_{j}$ then
\begin{equation}\label{*f11}
1\approx \frac{\olll(X'(j,\nk+\ooo_{j}))}
{\olll(\isss_{k_{j}}(\fff_{\nk})\sm
\isss_{k_{j-1}}(\fff_{n_{\kkks}}))}>\frac{1}2.
\end{equation}
Similarly to \eqref{2*f7a}
this estimate does not depend on the choice of $n_{s}$.
Denote by $X(1/2,3/4,j)$ the set of those $x$ for which
there exists $t\in\Z$ such that
$x\in[ (t+\frac{1}2)2^{\nk+k_{j-1}},
(t+\frac{3}4)2^{\nk+k_{j-1}}).$
We split $X(1/2,3/4,j)$ into two subsets depending on the parity
of $t$. If $t$ is even then $x\in X(1/2,3/4,j,e)$ and if
$t$ is odd then $x\in X(1/2,3/4,j,o)$.
Suppose $x,y\in X(1/2,3/4,j)$, $x<y$,
$x\in X(1/2,3/4,j,e)$, $y\in X(1/2,3/4,j,o)$ and $0<l=\lf y\rf -\lf x\rf$.
Then, as the reader can verify, $\lf x\rf\in X(1/2,3/4,j,e)$,
$\lf y\rf \in X(1/2,3/4,j,o)$, $y'\defeq x+l\in X(1/2,3/4,j,o)$
and
$x+2l\in [2t'\cdot 2^{\nk+k_{j-1}},(2t'+1)\cdot 2^{\nk+k_{j-1}})$
with a $t'\in\Z$.

Hence,
\begin{align}\label{2*f11}
&\text{if }x\in X'(j,\nk+k')\cap X(1/2,3/4,j,e),\
y\in X(1/2,3/4,j,o), \text{ and }\\
\nonumber &
0<l=\lf y\rf -\lf x\rf \leq 2^{\nk+k'}\leq 2^{\nk+\ooo_{j}}\text{ then}\\
& \nonumber x\in\isss_{k_{j}}(\fff_{\nk})\sm\isss_{k_{j-1}}(\fff_{\nk}),\
y'=x+l\in\isss_{k_{j-1}}(\fff_{\nk}),\text{ and }\\
&\nonumber
x+2l\in \isss_{k_{j}}(\fff_{\nk})\sm\isss_{k-1}(\fff_{\nk}).
\end{align}
Moreover, from $x\in  [(t+\frac{1}2)2^{\nk+k_{j-1}},
(t+\frac{3}4)2^{\nk+k_{j-1}})$ with $t$ even it follows that
\begin{equation}\label{*81*}
[x-2^{n_{s}},x+2^{n_{s}}+1)\sse [t2^{n_{s}+k_{j-1}},(t+1)2^{n_{s}+k_{j-1}})
\sse \isss_{k_{j}}(\fff_{\nk})\sm\isss_{k_{j-1}}(\fff_{\nk}).
\end{equation}

\subsection{SUBSTEP 2d: Estimate of the measure of those
points where \eqref{2*f10} holds in one $2^{n_{s}}$
grid interval, definition of the sets
$\cax_{t_{0},t_{1},k'}$}

Suppose
\begin{equation}\label{*f12}
[t_{0}2^{\nk},(t_{0}+1)2^{\nk})\sse X'(j,\nk+\ooo_{j})\cap
X(1/2,3/4,j,e).
\end{equation}
By \eqref{*81*} used with $x=t_{0}2^{n_{s}}$
we have
\begin{equation}\label{*82*}
[(t_{0}-1)2^{n_{s}},(t_{0}+1)2^{n_{s}}+1)\sse
\isss_{k_{j}}(\fff_{\nk})\sm\isss_{k_{j-1}}(\fff_{\nk}).
\end{equation}
We want to obtain an estimate of the measure of those $x$'s
for which \eqref{2*f10} and hence \eqref{*f10} holds for a suitable $l$.

By \eqref{*f12} there exists $t_{0}'\in \Z$ such that
\begin{equation}\label{4*f7cf12}
[t_{0}2^{\nk}, (t_{0}+1)2^{\nk})\sse
[(2t_{0}'+\frac{1}2)2^{\nk+k_{j-1}},
(2t_{0}'+\frac{3}4)2^{\nk+k_{j-1}}).
\end{equation}
Suppose
\begin{align}\label{**f7cf12}
&y\in X(1/2,3/4,j,o),\ 0<l=\lf y\rf -\lf x\rf \leq
2^{\nk+k'}\leq 2^{\nk+\ooo_{j}}\text{ and }\\
&\nonumber
y\in\supp(\fff_{\nk}).
\end{align}

By \eqref{**f7cf12},
 $\fff_{\nk}(y)=2^{\nk}$ and by \eqref{*11ns}
\begin{equation}\label{4*f7cf12b}
\exists\, t_{1}\in \Z\text{ such that }
y\in [t_{1}2^{\nk}, t_{1}2^{\nk}+1).
\end{equation}
Moreover, $l=t_{1}2^{n_{s}}-\lf x \rf$
and $y'=x+l\in [t_{1}2^{n_{s}}, t_{1}2^{n_{s}}+1)$
as well, which implies $\fff_{n_{s}}(x+l)=\fff_{n_{s}}(y')=2^{n_{s}}$.
Denote by
\begin{equation}\label{3*f7cf13}
\Ics j(x)=[(\tcs j(x)-1)2^{\nk+k_{j}},
\tcs j(x)2^{\nk+k_{j}})
\end{equation}
the component of $\isss_{k_{j}}(\fff_{\nk})$ containing $x$.

By the definition of $X'(j,\nk+\ooo_{j})$,
$y'=x+l$ and $x+2l$ are in $\Ics j(x)$, that is, $x,\, y'=x+l,$
and $x+2l$ belong to the same component of $\isss_{k_{j}}(\fff_{\nk})$.
Moreover, by \eqref{2*f11}, $y'=x+l\in\isss_{k_{j-1}}(\fff_{\nk})$
and $x+2l\in\isss_{k_{j}}(\fff_{\nk})\sm \isss_{k_{j-1}}(\fff_{\nk})$.
We introduce the notation
$$\Ics{j-1,1/2,3/4}(y)=[(\tcs{j-1,1/2,3/4}(y)+\frac{1}2)2^{\nk+k_{j-1}},
(\tcs{j-1,1/2,3/4}(y)+\frac{3}4)2^{\nk+k_{j-1}})
$$
for the interval containing $y\in X(1/2,3/4,j)$.
Then $y'\in \Ics{j-1,1/2,3/4}(y)$ holds as well.
Recall that we still suppose $j\geq 2.$
By our construction, see Section \ref{subsubstep1a} and
the second paragraph of \ref{subsubstep2a} we have
\begin{align}\label{3*f7df13}
&\text{if }I^{*}=[(t^{*}+\frac{1}2)2^{\nk+k_{j-1}},
(t^{*}+\frac{3}4)2^{\nk+k_{j-1}})\text{ is a component of}\\
& \nonumber
\isss_{k_{j-1}}(\fff_{\nk})\cap X(1/2,3/4,j)\text{ then}\\
\label{*f7df13} &
\lll(\isss_{0}(\fff_{\nk})\cap I^{*})=
\oot{j-1}\lll(I^{*})=\oot{j-1}2^{\nk+k_{j-1}-2}.
\end{align}
This implies that there exist
\begin{equation}\label{**f7df14}
\oot{j-1}2^{k_{j-1}-2}\text{ many }t_{1}\text{'s such that }
[t_{1}2^{\nk},t_{1}2^{\nk}+1)\sse\supp(\fff_{\nk})\cap I^{*}.
\end{equation}
We denote the set of these $t_{1}$'s by $\cat_{1}(I^{*})$. In the
interval $[t_{0}2^{\nk}+2^{\nk+k'-1}, t_{0}2^{\nk}+2^{\nk+k'})$
there exist $2^{k'-k_{j-1}-2}$ many $I^{*}$'s satisfying
\eqref{3*f7df13}. Denote the set of the corresponding $t^{*}$'s by
$\cat^{*}(t_{0}, k')$. Finally, denote by $\cat_{1}(t_{0},k')$ the
set of those $t_{1}$ which belong to a $\cat_{1}(I^{*})$ with
$I^{*}$ of the form in \eqref{3*f7df13} and
$t^{*}\in\cat^{*}(t_{0},k')$. Then
\begin{equation}\label{*f7ef14}
\#\cat_{1}(t_{0},k')=\oot{j-1}2^{k_{j-1}-2}\cdot 2^{k'-k_{j-1}-2}=
\oot{j}2^{k'-3}.
\end{equation}

Next, suppose $t_{1}\in \cat_{1}(t_{0},k')$ is fixed.
Then, from $y\in[t_{1}2^{\nk},t_{1}2^{\nk}+1)$ it follows that
$\fff_{\nk}(y)=2^{\nk}$.
For $x\in[t_{0}2^{\nk},(t_{0}+1)2^{\nk})$ set
\begin{equation}\label{3*f15}
\lxto=t_{1}2^{\nk}-\lf x \rf=\lf y \rf-\lf x \rf.
\end{equation}
Then $\fff_{n_{s}}(x+\lxto)=\fff_{n_{s}}(y')=\fff_{n_{s}}(y)=2^{n_{s}}.$
From $t_{1}\in \cat_{1}(t_{0},k')$ it follows that
\begin{equation}\label{*74*}
[t_{1} 2^{n_{s}},
t_{1}2^{n_{s}}+1)\sse I^{*}
\sse
[t_{0}2^{n_{s}}+2^{n_{s}+k'-1},t_{0}2^{n_{s}}+2^{n_{s}+k'})
\end{equation}
with a suitable $I^{*}$.
This implies
\begin{equation}\label{*67b}
0<t_{1}2^{n_{s}}-(t_{0}+1)2^{n_{s}}+1\leq l_{x,t_{1}}\leq t_{0}2^{n_{s}}+
2^{n_{s}+k'}-t_{0}2^{n_{s}}=2^{n_{s}+k'}.
\end{equation}
We have  $\fff_{\nk}(x+\lxto)=2^{\nk}$ and
\begin{equation}\label{*18**}
\frac{\fff_{\nk}(x+\lxto)\ppp_{k'}(x+2\lxto)}{\lxto}\geq 1
\end{equation}
holds if $x+2\lxto\in\supp(\ppp_{k'})$.

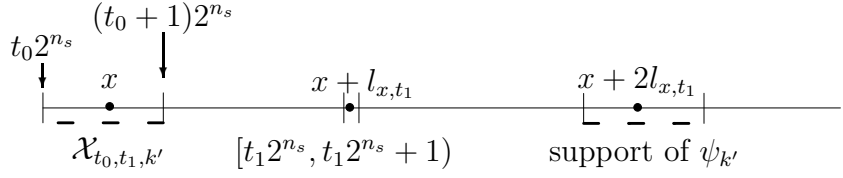
\placedrawing[h]{fgt5.lp}{The set $\cax_{t_{0},t_{1},k'}$}{figure5}

Denote by $\cax\sto{t_{1},k'}$ the set of those $x\in[t_{0}2^{\nk} ,
(t_{0}+1)2^{\nk} )$ for which $x+2\lxto\in\supp(\ppp_{k'})$.
By \eqref{*74*}
\begin{equation}\label{*62*lxt}
\text{for fixed }t_{0},\, t_{1}\text{ there can be at most one }k'
\text{ for which }t_{1}\in \cat_{1}(t_{0},k').
\end{equation}
Moreover,
\begin{equation}\label{*62*2}
\text{if }x\in\cax_{t_{0},t_{1},k'}\sse [t_{0}2^{n_{s}},(t_{0}+1)2^{n_{s}})
\text{ then}
\end{equation}
$$\text{\eqref{*18**} holds and }x+\lxto\in[t_{1}2^{n_{s}},t_{1}2^{n_{s}}+1).$$

Dividing \eqref{*67b} by $2^{n_{s}}$ and rearranging
we obtain that for $t_{1}\in \cat_{1}(t_{0},k')$
$$0< t_{1}-t_{0}\leq 2^{k'}+1-\frac{1}{2^{n_{s}}},$$
since $t_{1}$ and $t_{0}$ are integers
recalling \eqref{*40b*} and \eqref{2*f10}
we have
\begin{equation}\label{*70b*}
0< t_{1}-t_{0}\leq 2^{k'}<2^{\ooo_{j}}<2^{k_{M}}.
 \end{equation}

\subsection{SUBSTEP 2e: ``Periodicity and independence"
of the sets $\cax_{t_{0},t_{1},k'}$}

Observe that if $x\in[t_{0}2^{\nk},(t_{0}+1)2^{\nk})$
then
$\lf x\rf \in [t_{0}2^{n_{s}},(t_{0}+1)2^{n_{s}}-1]$
and by \eqref{3*f15},
$x+2\lxto\in [(2t_{1}-t_{0}-1)2^{\nk}+1,
(2t_{1}-t_{0})2^{\nk}+1)$, moreover by \eqref{*f9c} in this interval
\begin{equation}\label{2*f15}
\ppp_{k'}\text{ is ``periodic" by }
\pi(2t_{1}-t_{0},j).
\end{equation}
This implies that
\begin{align}\label{*f7ff15}
&\text{if }x,x+\pi(2t_{1}-t_{0},j)\in [t_{0}2^{\nk},(t_{0}+1)2^{\nk})
\text{ then }\\&
\nonumber x\in\cax\sto{t_{1},k'}\text{ iff }x+\pi(2t_{1}-t_{0},j)\in
\cax\sto{t_{1},k'}.
\end{align}

The ``periodic density"
of the measure of the support of
$\ppp_{2t_{1}-t_{0},j,k'}$
in
$[(2t_{1}-t_{0}-1)2^{\nk}+1,
(2t_{1}-t_{0})2^{\nk}+1)$
 equals $2^{-k'}\ointt\ppp_{2t_{1}-t_{0},j,k'}$.
So the measure of this support in an interval
of length $2^{n_{\kkks}}$
is approximately $(\ointt\ppp_{2t_{1}-t_{0},j,k'})2^{-k'}2^{n_{\kkks}}$
provided $2^{n_{s}}$ is much larger than $\pi(2t_{1}-t_{0},j)$.
Hence,
we have by \eqref{*f16} that
\begin{equation}\label{4*f16}
\frac{1}4 \left (\ointt\ppp_{
2t_{1}-t_{0},j,k'} \right )2^{-k'}2^{\nk}<
\left (\ointt\ppp_{
2t_{1}-t_{0},j,k'} \right )2^{-k'}2^{\nk}\approx
\lll(\cax\sto{t_{1},k'})<
\end{equation}
$$
2 \left (\ointt\ppp_{
2t_{1}-t_{0},j,k'} \right )2^{-k'}2^{\nk}.
$$
Suppose $t_{1}\in\cat_{1}(t_{0},k')$. Denote by
$\cak(t_{0},t_{1},k')$ the set of those $(t_{1}',k'')
\not=(t_{1},k')$
for which $t_{1}'\in\cat_{1}(t_{0},k'')$.
Observe that if $t_{1}=t_{1}'$ then from
\eqref{*62*lxt} it follows that for $k''\not=k'$
we have $t_{1}\not\in \cat_{1}(t_{0},k'')$.
Hence for $(t_{1}',k'')\in \cak(t_{0},t_{1},k')$
we have $t_{1}'\not=t_{1}$.
By \eqref{2*f9}, \eqref{*70b*} and \eqref{*f7ff15} if
$(t_{1}',k'')\in\cak(t_{0},t_{1},k')$ then the ``periods" of
$\cax\sto{t_{1},k'}$ and $\cax\sto{t_{1}',k''}$ are relatively
prime and these sets are ``independent" in the sense that,
using \eqref{4*f16} and Lemma \ref{*indep} for large $n_{s}$ as well,
we have
\begin{align}\label{*f7gf17}
&\lll(\cax\sto{t_{1},k'}\cap \cax\sto{t_{1}',k''})<
\lll(\cax\sto{t_{1},k'})\cdot 2\cdot \frac{\lll(\cax\sto{t_{1}',k''})}
{2^{\nk}}<\\&\nonumber
4\cdot 2^{-k''}(\ointt\ppp_{2t_{1}'-t_{0},j,k''})\lll(\cax\sto{t_{1},k'}).
\end{align}
Using  \eqref{3*f16}, \eqref{*f7ef14}, \eqref{*62*lxt} and \eqref{*f7gf17} we obtain
\begin{align}\label{2*f7gf17}
&\lll(\cax\sto{t_{1},k'}\sm \bigcup_{(t_{1}',k'')\in\cak(t_0,t_{1},k')}\cax
\sto{t_{1}',k''})>\\&\nonumber
\left (1-\sum_{k''=\aaa_{j}}^{\ooo_{j}}4\cdot 2^{-k''}\left
(\max_{t}\ointt\ppp_{t,j,k''}\right)
\oot j2^{k''-3}
\right)\lll(\cax\sto{t_{1},k'})>\frac{3}4\lll(\cax\sto{t_{1},k'}).
\end{align}

\subsection{SUBSTEP 2f: Estimates of the measure
of points where the maximal operator is large for the
$(s-1)$-functions, the sets $\cax_{t_{0},s-1}$}

We recall that by \eqref{2*f11}, \eqref{*f12} and \eqref{*82*}
$$
[t_{0}2^{\nk} ,
(t_{0}+1)2^{\nk}+1)\sse \isss_{k_{j}}(\fff_{\nk})\sm\isss_{k_{j-1}}(
\fff_{\nk}).
$$
By \eqref{*73oo*}, \eqref{*44b} and \eqref{*73xx*}
\begin{align}\label{*f19}
&\text{on }
[t_{0}2^{\nk}+1,
(t_{0}+1)2^{\nk}+1 )\text{ the functions }\\
&\nonumber \fff_{i}, \ppp_{k}, i=1,...,n_{j,\kkks-1},
k=\aaa_{j},...,\ooo_{j} \text{ are the restrictions of an}\\&
\nonumber \aaa_{j}-\capp(t_{0}+1,j)-(\kkks-1)-\text{family
 periodic by }\pi(t_{0}+1,j).
\end{align}
Set
\begin{align}\nonumber
\cax\sto{\kkks-1}=& \{ x\in
[t_{0}2^{\nk} ,
(t_{0}+1)2^{\nk}):\\&\label{**19*}
\max_{k\in[\aaa_{j},\ooo_{j}]
}\max _{l\leq\pi^{*}}\frac{\sum_{i=1}^{n_{j,\kkks-1}}
\fff_{i}(x+l)\ppp_{k}(x+2l)}{l}\geq 1
  \}.
\end{align}
These are the ``old" points in $[t_{0}2^{n_{s}}, (t_{0}+1)2^{n_{s}})$
where the inequality from \eqref{**19*} holds. We will need those ``new"
points from the sets $\cax_{t_{0},t_{1},k'}$ which do not belong
to these ``old" sets
$\cax_{t_{0},s-1}$.
By \eqref{*f19}
\begin{align}\label{2*f19}
&\text{if }x,x+\pi(t_{0}+1,j)\in [t_{0}2^{\nk}+1,
(t_{0}+1)2^{\nk}-2\pi^{*})\text{ then }\\&\nonumber
x\in\cax\sto{\kkks-1}\text{ iff }x+\pi(t_{0}+1,j)\in\cax\sto{\kkks-1}.
\end{align}
In \eqref{*48b*} we supposed that $
\pi(t_{0}+1,j)\leq \pi^{*}\ll 2^{\nk}$, so for most part of
$[t_{0}2^{\nk},(t_{0}+1)2^{\nk})$ the set $\cax\sto{\kkks-1}$
is ``periodic" by $\pi(t_{0}+1,j)$.

By our assumptions  for our
$\aaa_{j}-\capp(t_{0}+1,j)-(\kkks-1)$-family $\olll \{ .  \}\leq
\eee<\eee_{0}=\frac{1}{10}$ holds in \eqref{*f1} therefore, for
sufficiently large $n_{s}$
\begin{equation}\label{*65a}
\lll(\cax\sto{\kkks-1})=\lll(\cax\sto{\kkks-1}\cap [t_{0}2^{\nk},(t_{0}+1)2^{\nk}))<\frac{1}8
2^{\nk}.
\end{equation}

\subsection{SUBSTEP 2g: Independence estimates of ``old"
and ``new" sets, that is, of $\cax_{t_{0},s-1}$ and $\cax_{t_{0},t_{1},k'}$.
Estimates of the ``new contribution" set $\cax'_{t_{0}}$}

By \eqref{2*f19}, $\cax\sto{\kkks-1}$ is ``periodic"
by $\pi(t_{0}+1,j)$  while by \eqref{*f7ff15}, $\cax\sto{t_{1},k'}$
is ``periodic" by $\pi(2t_{1}-t_{0},j).$
By \eqref{*40b*} and \eqref{*70b*}
$$0<2t_{1}-t_{0}-(t_{0}+1)=2t_{1}-2t_{0}-1<2^{k'+1}\leq 2^{\ooo_{j}+1}
\ll 2^{k_{M}+1}.$$
Hence by \eqref{2*f9}, $\pi(t_{0}+1,j)$ and $\pi(2t_{1}-t_{0},j)$ are
relatively prime and the sets $\cax\sto{t_{1},k'}$ are sufficiently
``independent" of $\cax\sto{\kkks-1}$.
By \eqref{*65a} and Lemma \ref{*indep}
for $2^{\nk}\gg \pi(t_{0}+1,j)$
we have
\begin{equation}\label{*675}
\lll(\cax\sto{t_{1},k'}\cap \cax\sto{\kkks-1})<\frac{1}4\lll
(\cax\sto{t_{1},k'}).
\end{equation}
Using \eqref{4*f16}, \eqref{2*f7gf17}
 and  \eqref{*675} we infer
\begin{align}\label{*f21}
&\lll\left (\cax\sto{t_{1},k'}\sm \left (
\cax\sto{\kkks-1}\cup
\bigcup_{(t_{1}',k'')\in\cak(t_0,t_{1},k')}\cax
\sto{t_{1}',k''}\right)\right)>\\&\nonumber
\frac{1}2\lll(\cax\sto{t_{1},k'})>\frac{1}8 \left (\min_{t}
\ointt\ppp_{t,j,k'}\right)
2^{-k'}2^{\nk},
\end{align}
this estimates the ``new contribution" of $\cax\sto{t_{1},k'}$
to the set of points where \eqref{2*f10} holds.
Set
\begin{equation}\label{*83plus*}
\cax'\sto{k'}=\bigcup_{t_{1}\in\cat_{1}(t_{0},k')}
\left (\cax\sto{t_{1},k'}\sm\left (\bigcup_{(t_{1}',k'')
\in\cak(t_0,t_{1},k')}\cax\sto{t_{1}',k''}\cup \cax\sto{\kkks-1}\right)\right).
\end{equation}
From $k'\not=k''$ and $x\in \cax'\sto{k'}\cap \cax'\sto{k''}$
it follows that $x\in \cax_{t_{0},t_{1},k'}$ with $t_{1}\in \cat_{1}(t_{0},k')$
and $x\in \cax_{t_{0},t_{1}',k''}$ with $t_{1}'\in \cat_{1}(t_{0},k'')$,
that is $(t_{1}',k'')\in \cak(t_{0},t_{1},k')$,
but this contradicts \eqref{*83plus*}.
Hence
\begin{equation}\label{*110b*}
\cax'\sto{k'}\cap \cax'\sto{k''}=\ess\text{ for }k'\not=k''.
\end{equation}

By \eqref{*f7ef14} and \eqref{*f21}
\begin{equation}\label{*f22}
\lll(\cax'\sto{k'})>\bigg(\frac{1}8\bigg
(\min_{t}\ointt\ppp_{t,j,k'}\bigg)2^{-k'}2^{\nk}\bigg)
\oot{j}2^{k'-3}=
\end{equation}
$$\frac{1}{64}\bigg(\min_{t}\ointt\ppp_{t,j,k'}\bigg)2^{\nk}\oot{j}.$$

Set
\begin{equation}\label{*111*}
\cax'_{t_{0}}=\cup_{k'=\aaa_{j}}^{\ooo_{j}}\cax'\sto{k'}.
\end{equation}
Then by $\rrr>1/2$, \eqref{3*f16}, \eqref{*110b*} and \eqref{*f22}
\begin{align}\label{2*f22}
&\lll(\cax'_{t_{0}})>\sum_{k'=\aaa_{j}}^{\ooo_{j}}\frac{1}{64}
\left (\min_{t}\ointt\ppp_{t,j,k'}\right)2^{\nk}\oot{j}=\\&\nonumber
\frac{1}{64}\left (\sum_{k'=\aaa_{j}}^{\ooo_{j}} \min_{t}
\ointt\ppp_{t,j,k'}\right)
2^{\nk}\oot{j}>\frac{1}{128}\eee2^{\nk}\oot{j}.
\end{align}
Hence in each interval $I(t_{0})\sse[t_{0}2^{\nk},(t_{0}+1)2^{\nk})$
satisfying \eqref{*f12} we have found a set $\cax'_{t_{0}}$
such that
\begin{equation}\label{*f7hf23}
\lll(\cax'_{t_{0}})>\frac{\eee}{128}\oot{j}\lll(I(t_{0}))
\end{equation}
and
for $x\in\cax'_{t_{0}}$
\begin{equation}\label{3*f7hf23}
\exists l\leq 2^{n_{s}+\ooo_{j}}
\text{ such that }
\frac{\fff_{\nk}(x+l)\sum_{k'=\aaa_{j}}^{\ooo_{j}}\ppp_{k'}(x+2l)}
{l}\geq 1.
\end{equation}
Moreover,
\begin{equation}\label{*87plus*}
\cax'_{t_{0}}\cap\cax\sto{\kkks-1}=\ess.\end{equation}

\subsection{SUBSTEP 2h: Conclusion}

Denote by $\cai(j)$ the union of all intervals in $\isss_{k_{j}}
(\fff_{\nk})\sm
\isss_{k_{j-1}}(\fff_{\nk})$ which satisfy \eqref{*f12}.
Then by \eqref{*9ns*}, \eqref{*10ns*}, \eqref{*f11}, and the definition of
$X(1/2,3/4,j,e)$ we have
\begin{equation}\label{2*f7hf24}
\olll(\cai(j))> \frac{1}{16}\olll(\isss_{k_{j}}(\fff_{\nk}))=
\frac{1}{16}2^{-M+j}.
\end{equation}
Denote by $\ocax(j)$ the set of those $x\in\isss_{k_{j}}(\fff_{\nk})
\sm
\big(\isss_{k_{j-1}}(\fff_{\nk})\cup\bigcup_{t_{0}}\cax\sto{\kkks-1}\big)$
for which there exists $l$ such that \eqref{3*f7hf23} holds.
Then \eqref{*f7hf23} and \eqref{2*f7hf24} imply that
\begin{equation}\label{*f7if24}
\olll(\ocax(j))>\frac{\eee}{128\cdot 16}\oot{j}2^{-M+j}=
\frac{\eee}{2048}2^{-M}.
\end{equation}
When $j\not=j'$
the sets $\ocax(j)\sse \isss_{k_{j}}(\fff_{\nk})\sm\isss_{k_{j-1}}
(\fff_{\nk})$ and
$\ocax(j')\sse \isss_{k_{j'}}(\fff_{\nk})\sm\isss_{k_{j'-1}}
(\fff_{\nk})$ are disjoint from each other and from the
sets $\cup_{t_{0}}\cax'\sto{\kkks-1}$.
Using the estimates \eqref{*f7if24} for $j=2,...,M$
and \eqref{*1043b}
we obtain that \eqref{*f1} holds with $s$ as well.
Inequality \eqref{**f1} follows from \eqref{*13ns*}, \eqref{*3rel},
and \eqref{**3rel}.
Inequality \eqref{***f1} follows from \eqref{*4rel} and \eqref{**4rel}.
\end{proof}

\subsection{\bf How to derive Theorem \ref{ith1} from Theorem \ref{f2}}
\label{derith1}
   \medskip

\medskip

   We define the dynamical system $(Y_s, \mathcal{B}_s,
   \overline{\lambda}_s, T_{p_s})$ where
   $Y_s= [0, p_s),$ $T_{p_s}: x\rightarrow x+ 1$ $\text{mod}\,  p_s$,
   $\mathcal{B}_s$ is the $\sigma$-field of measurable subsets on
   $Y_s$ and $\olll_{s}$ is the normalized Lebesgue measure on
$[0,p_{s}).$

By Theorem \ref{f2} we have
   $$\overline{\lambda}_s\bigg\{x: \max_{l\leq
   p_s}\sum_{k=\alpha}^{\omega}\sum_{i=1}^{n_s}
   \frac{\phi_i(T^l_{p_s}(x))\psi_k(T^{2l}_{p_s}(x))}{l}\geq 1 \bigg\}>
   \min\big\{\epsilon,\frac{s\epsilon(M-1)2^{-M}}{2048}\big\}.$$

 Dividing by the norms of the sums of the non-negative functions
$\psi_k$ and $\phi_i$ we derive the estimate
 $$\overline{\lambda}_s\bigg\{x: \max_{l\leq
 p_s}\sum_{k=\alpha}^{\omega}\sum_{i=1}^{n_s}
   \frac{\phi_i(T^l_{p_s}(x))\psi_k(T^{2l}_{p_s}(x))}{l\|\sum_{i=1}^{n_s}\phi_i\|_1\|
\sum_{k=\alpha}^{\omega}\psi_k\|_1}\geq
   \frac{1}{\|\sum_{i=1}^{n_s}\phi_i\|_1\|\sum_{k=\alpha}^{\omega}\psi_k\|_1}\bigg\}>$$
 $$  \min\big\{\epsilon,\frac{s\epsilon(M-1)2^{-M}}{2048}\big\}.$$

  For $0<\epsilon <\epsilon_0$
fixed, one can pick $s$ and $M$ such that
  $\frac{s(M-1)2^{-M}}{2048}\geq 1$.
Actually we will choose $s$ to be equal to
   $\dis \lfloor\frac{2048\cdot 2^M}{M-1}\rfloor + 1.$
We can observe that with this choice for $M>2$ we have
  $\dis \lfloor\frac{2048\cdot 2^M}{M-1}\rfloor + 1 \leq
  2\frac{2048\cdot 2^M}{M-1}.$
  Using the estimates \eqref{**f1} and \eqref{***f1}
 we have
   $$\overline{\lambda}_s\bigg\{x: \max_{l\leq
 p_s}\sum_{k=\alpha}^{\omega}\sum_{i=1}^{n_s}
   \frac{\phi_i(T^l_{p_s}(x))\psi_k(T^{2l}_{p_s}(x))}
{l\|\sum_{i=1}^{n_s}\phi_i\|_1\|\sum_{k=\alpha}^{\omega}\psi_k\|_1}\geq
   \frac{1}{\epsilon s2^{-M+1}}\bigg\}> \epsilon.$$

   Therefore, if we set $\dis F_s\defeq
\frac{\sum_{i=1}^{n_s}\phi_i}{\|\sum_{i=1}^{n_s}\phi_i\|_1}$ and
   $\dis G_s\defeq \frac{\sum_{k= \alpha}^{\omega}
   \psi_k}{\|\sum_{k=\alpha}^{\omega}\psi_k\|_1}$ then we have
   $$\overline{\lambda}_s\bigg\{x: \max_{l\leq
 p_s} \frac{F_s(T^{l}_{p_s}(x))G_s(T^{2l}_{p_s}(x))}{l}\geq
 \frac{1}{\epsilon s2^{-M+1}}\bigg\}> \epsilon.$$
 As we have assumed that $\dis s =
\lfloor\frac{2048\cdot 2^M}{M-1}\rfloor + 1,$
 the last inequality implies
   \begin{equation}\label{idris1}
   \overline{\lambda}_s\bigg\{x: \max_{l\leq
 p_s} \frac{F_s(T^{l}_{p_s}(x))G_s(T^{2l}_{p_s}(x))}{l}\geq \frac
 {(M-1)}{4\cdot 2048\epsilon}\bigg\}>\epsilon.
 \end{equation}
 To simplify further the notation we set $C=
 \frac{1}{8192\epsilon}.$ This constant depends only on
 $\epsilon$,
a fixed positive number.
We introduce different values of $M$ with the sequence $M_j = (j+1)^5.$
 To these values
the subsequence $s_j$ of natural numbers
is automatically associated,  where
 $\dis s_j = \dis \lfloor\frac{2048\cdot 2^{(j+1)^5}}{(j+1)^5-1}
\rfloor +
 1.$
 We define now the measure preserving system
 $\dis \prod_{j=1}^{\infty}(Y_{s_j}, \mathcal{B}_{s_j},
   \overline{\lambda}_{s_j}, T_{p_{s_j}})$ and the functions
   $ \dis W(y_1, y_2, ...,y_n,...) = \sum_{j=1}^{\infty}
   \sqrt{\beta_j}F_{s_j}(y_j),
$ \\
   $\dis Q(y_1, y_2, ...,y_n,...) = \sum_{j=1}^{\infty}
   \sqrt{\beta_j} G_{s_j}(y_j)$ where the constants $\beta_j$ will be
   specified a bit later. The space is $\dis Y=\prod_{j=1}^{\infty}Y_{s_j}$,
   and the measure
preserving transformation acting on $Y$ is the map $U$
   defined as
   $U(y_1, y_2,..., y_n,...) = (T_{p_{s_1}}(y_1),
   T_{p_{s_2}}(y_2),...,T_{p_{s_n}}(y_n),...).$ The invariant measure under $U,$ is denoted by $\mu,$ and it is
   the product of the measures $\overline{\lambda}_{s_j}.$
   We have for any positive number $z,$
   $$
   \mu\bigg\{y\in Y: \sup_l \frac{W(U^l(y))Q(U^{2l}(y))}{l}\geq
   z\bigg\}\geq
$$ $$\overline{\lambda}_{s_j}\bigg\{y_j\in Y_{s_j}:
   \beta_j\sup_l\frac{F_{s_j}(T^l_{p_{s_j}}(y_j))G_{s_j}(T^{2l}_{p_{s_j}}(y_j))}{l}\geq z\bigg\}.
   $$
 We can select now the terms
 $\beta_j$. We put $\beta_j = (\frac{\pi^2}{6}-1)^{-2}(j+1)^{-4}.$
 Observe that this guarantees that
the functions $W$ and $Q$
 have norm $1$ in $L^1(\mu).$
  We also have $\dis \lim_j(M_{j}-1)\beta_j
= \infty.$
   With these conditions we can conclude without difficulty. For each
   $j$ we have
 $$
   \mu\bigg\{y\in Y: \sup_l \frac{W(U^l(y))Q(U^{2l}(y))}{l}\geq
   C\beta_j(M_j-1)
   \bigg\}\geq $$
   $$\overline{\lambda}_{s_j}\bigg\{y_{j}\in Y_{s_j}:
\sup_l
\frac{F_{s_j}(T^l_{p_{s_j}}(y_j))G_{s_j}(T^{2l}_{p_{s_j}}(y_j))}{l}\geq
 C (M_j-1)\bigg\}>\epsilon.
 $$
   Therefore, we have found a complete non-atomic measure space, a
   measure preserving system and functions $W$ and $Q$ with $L^1$
   norm $1$ such that
   $$\mu\bigg\{y\in Y: \sup_l \frac{W(U^l(y))Q(U^{2l}(y))}{l}= \infty
   \bigg\}>\epsilon.
   $$
   Thus Theorem \ref{ith1}
holds for this measure preserving system. By using the
  disintegration of $\mu$ into ergodic components i.e. $d\mu =
  d\mu_c dc$ and for a.e. c  and $U$ with respect to $\mu_c$ is
  ergodic we have
  $$\int \mu_c\bigg\{y\in Y: \sup_l \frac{W(U^l(y))Q(U^{2l}(y))}{l}= \infty
   \bigg\}dc> \epsilon.$$ Therefore there exists at least one $c^*$
   such that
  $$\mu_{c^*}\bigg\{y\in Y: \sup_l \frac{W(U^l(y))Q(U^{2l}(y))}{l}= \infty
   \bigg\}> \epsilon.$$
 To duplicate this for any ergodic measure preserving system
 we can
 use Halmos's result \cite{Halmos} on the density of the conjugates of an ergodic measure
 preserving system for the weak topology.  We have
 for any ergodic measure preserving system $(B, \mathcal{F}, m,
 D)$ on a complete non-atomic finite measure space, for any $j$
 \[
 \begin{aligned}
 &\sup_{\|f\|_1=1, \,\|g\|_1=1}m\bigg\{b:
 \sup_l\frac{f(D^lb)g(D^{2l}b)}{l}> M'_j\bigg\} \\
 &\geq \sup_{\|W\|_1=1,\,\|Q\|_1= 1}
 \mu_{c^*}\bigg\{y\in Y: \sup_l \frac{W(U^l y)Q(U^{2l}
  y )}{l}> M'_j
   \bigg\}>\epsilon.
   \end{aligned}
   \]
 The details are
 omitted here. Similar details are given in the appendix of \cite{Demeter}
 in the case of the averages associated with the triple a.e. recurrence
 and these arguments apply to the situation considered here.


 \end{document}